\documentclass[10pt]{article}
\usepackage{amsmath}
\usepackage{latexsym}
\usepackage{amsfonts}
\usepackage{amssymb}
\usepackage{amsthm}
\usepackage{amsbsy}
\usepackage{graphicx}
\usepackage{bbm,a4wide}
\usepackage{hyperref}
\usepackage{placeins}
\usepackage{ulem} 
\usepackage{color}
\usepackage{epstopdf}

\newcommand{\R}{\mathbb{R}}
\newcommand{\N}{\mathbb{N}}

\newcommand{\Per}{P_{\mathbb{H}}}

\newcommand{\Div}{\mathrm{div}}
\newcommand{\Var}{\mathrm{Var}}

\newcommand{\eps}{\varepsilon}

\newcommand{\opG}{\circ} 

\newtheorem{theorem}{Theorem}[section]
\newtheorem{lemma}[theorem]{Lemma}
\newtheorem{proposition}[theorem]{Proposition}
\newtheorem{definition}[theorem]{Definition}

\newtheorem{remark}[theorem]{Remark}

\begin{document}
\title{A semi-discrete scheme for horizontal mean curvature flow \\ in the Heisenberg group}

\author{
Annalisa Cesaroni\thanks{Dipartimento di Matematica ``Tullio Levi-Civita'', Universit\`{a} di Padova, via Trieste 63, 35131 Padova, Italy, \url{ annalisa.cesaroni@unipd.it}},
Matteo Novaga\thanks{Dipartimento di Matematica, Universit\`{a} di Pisa, Largo Bruno Pontecorvo 5, 56127 Pisa, Italy, \url{matteo.novaga@unipi.it}}~ and 
Paola Pozzi\thanks{Fakult\"at f\"ur Mathematik, Universit\"at Duisburg-Essen,
Thea-Leymann-Stra\ss e 9,
45127 Essen, Germany, \url{paola.pozzi@uni-due.de}}
}

\date{}

\maketitle

\begin{abstract}
The horizontal mean curvature flow is a geometric evolution associated with the horizontal perimeter in sub-Riemannian spaces. 
In this paper, we introduce a minimizing movement scheme for this flow in the first Heisenberg group $\mathbb{H}$, and prove the existence of a flat flow which is locally $1/2$-H\"older continuous in time, with respect to the $L^1$-distance. The construction applies to arbitrary bounded sets of finite horizontal perimeter and does not require smoothness or symmetry of the initial datum.
\end{abstract}

\tableofcontents

\section{Introduction}

The Heisenberg group is the simplest noncommutative Carnot group and provides a basic setting for the study of geometric variational problems in sub-Riemannian geometry \cite{CapognaDanielli, mont, Monti}. Its geometry is determined by a two-dimensional, nonintegrable horizontal distribution and by the associated Carnot--Carath\'eodory distance. The horizontal perimeter measures the variation of characteristic functions along the horizontal directions. The theory of sets of finite horizontal perimeter, including rectifiability and compactness, supplies the measure-theoretic framework for these problems \cite{serra, CapognaDanielli}. The group $\mathbb{H}$ has homogeneous dimension $4$: the volume scales with degree $4$, whereas the horizontal perimeter scales with degree $3$.

For a smooth surface $\Sigma_t=\partial E_t$, horizontal mean curvature flow prescribes, at non-characteristic points, the velocity
\[
\partial_t p(t)=-H_{\mathbb{H}}(p(t),E_t)\nu_{\mathbb{H}}(p(t)),
\]
where $\nu_{\mathbb{H}}$ is the horizontal outer unit normal and $H_{\mathbb{H}}=\Div_{\mathbb{H}}\nu_{\mathbb{H}}$ is its horizontal divergence. A distinctive difficulty is the presence of characteristic points, at which the tangent plane coincides with the horizontal distribution. At such points the horizontal normal is not defined, even when the surface is smooth in the Euclidean sense. Thus, besides the possible formation of singularities during the evolution, the geometry itself introduces a degeneracy in the classical motion law.

Related difficulties already arise in the stationary theory. Danielli, Garofalo, Nhieu and Pauls \cite{dgar} show that stability restores a Bernstein-type rigidity: the only stable $C^2$ entire $H$-minimal graphs with empty characteristic locus are vertical planes. For the isoperimetric problem, existence of minimizers in Carnot groups was established by Leonardi and Rigot \cite{lr}. 
Pansu balls are the natural candidate minimizers in $\mathbb{H}$, and their role is supported by classification results under axial symmetry, Euclidean convexity, or $C^2$ regularity assumptions \cite{m, mr, rr}. In particular, Monti and Rickly \cite{mr} characterize convex isoperimetric sets without imposing additional smoothness or symmetry. Calibration methods also yield isoperimetric comparisons with Pansu balls \cite{r}. These results provide useful geometric models and variational tools for the evolution problem considered here.

Several weak approaches to horizontal mean curvature flow have been developed. In the level-set formulation, the evolving sets are represented as sublevel sets of a function satisfying, away from points where its horizontal gradient vanishes,
\[
u_t=|\nabla_{\mathbb{H}}u|\,
\Div_{\mathbb{H}}\, \frac{\nabla_{\mathbb{H}}u}{|\nabla_{\mathbb{H}}u|}\,. 
\] Capogna and Citti \cite{capcitti2} extended the viscosity approach to Carnot groups, proving existence and comparison results under suitable restrictions. Ferrari, Liu and Manfredi \cite{ferrari} established existence, uniqueness and stability for axisymmetric viscosity solutions in $\mathbb{H}$, and obtained an explicit evolution starting from a Kor\'anyi sphere. Baspinar and Citti \cite{bc} proved comparison and uniqueness results in step-two Carnot groups and in the roto-translation group by relating viscosity solutions to Riemannian approximations; their comparison theorem is formulated for continuous solutions which are constant outside a compact set.  A paper   by Bombardieri, Fogagnolo and Franceschi \cite{gaia}  exploring  qualitative properties of  level set viscosity solutions with  mean-convex initial data is in preparation.  Stochastic control representations provide another approach to the level-set equation \cite{dirr}, with further developments in \cite{ddg, grande} concerning optimal controls and approximating flows. More recently, Citti, Dirr, Dragoni and Grande \cite{cddg} derived horizontal mean curvature motion through a formal asymptotic expansion of a nonlocal mean-field equation under Heisenberg dilations. Their analysis at non-characteristic points and numerical experiments provide a further approximation of the flow, motivated also by models of visual processing.

Our approach in this paper is variational. We adapt the time-discrete construction introduced in the Euclidean setting by Almgren, Taylor and Wang \cite{ATW} and Luckhaus and Sturzenhecker \cite{lu}. This construction has been developed for a variety of geometric evolutions, including forced, anisotropic and inhomogeneous mean curvature flows; see, for instance, \cite{BCK, ChamDeG}. Given a time step $h>0$ and a bounded set $E$ of finite horizontal perimeter, we consider minimizers of
\[
J_h(F):=P_{\mathbb{H}}(F)+\frac{1}{h}\int_F d_E(p)\,dp,
\]
where $d_E$ denotes the signed Carnot--Carath\'eodory distance, negative in $E$. Up to a constant independent of $F$, this functional can be written as
\[
P_{\mathbb{H}}(F)+\frac{1}{h}\int_{E\triangle F}|d_E(p)|\,dp.
\]
The second term penalizes the displacement from the preceding set, while minimization decreases the horizontal perimeter. Iterating the scheme, with a consistent choice of maximal or minimal minimizers, defines a piecewise constant evolution $E_h(t)$.

The use of the Carnot--Carath\'eodory distance is closely related to the geometry of the prescribed motion. Arcozzi and Ferrari \cite{af1, af2} studied metric normals and the regularity and horizontal Hessian of the signed distance to a smooth surface. In particular, for a $C^3$ boundary, at non-characteristic points one has
\[
\nabla_{\mathbb{H}}d_E=\nu_{\mathbb{H}},
\qquad
\Div_{\mathbb{H}}\nabla_{\mathbb{H}}d_E=H_{\mathbb{H}}.
\]
These identities motivate the relation between the discrete variational problem and the horizontal mean curvature flow. They also indicate why consistency is delicate near characteristic points, where a regular tubular description by the signed distance is unavailable.

An important ingredient of our analysis is a class $\mathcal{B}$ of bounded $C^2$, axially symmetric sets which are strictly starshaped with respect to the group of dilations, horizontally mean convex, and have bounded horizontal mean curvature; see Definition~\ref{definitionb}. Extending the horizontal normal along the foliation generated by dilations, we construct a calibration and prove that these sets are outward minimizing. More precisely, if $H_{\mathbb{H}}\geq\delta\geq0$ on the non-characteristic part of $\partial B$, $B \in \mathcal{B}$, Proposition~\ref{convexset} gives
\[
P_{\mathbb{H}}(F)\geq P_{\mathbb{H}}(B)
+\frac{\delta}{\lambda_F}|F\setminus B|,
\qquad B\subseteq F,
\]
where $\lambda_F$ is the smallest dilation parameter needed to contain $F$, with $\lambda_F=+\infty$ for unbounded $F$. The class includes Kor\'anyi balls, with $\delta=0$, and Pansu balls of parameter $R$, with $\delta=2/R$. In particular, Pansu balls give uniform confinement of the discrete evolutions. We also establish order preservation for the extremal minimizers and preservation of outward minimality under the hypotheses of Proposition~\ref{outward}. For strictly mean convex members of $\mathcal{B}$, Proposition~\ref{speedpansu} gives a quantitative one-step contraction estimate. The symmetry assumptions enter these geometric constructions, but are not imposed on the sets evolved by the scheme.

We then prove estimates for general discrete minimizers $F\in T_h(E)$ of the functional $J_{h}$. The displacement bound
\[
\sup_{E\triangle F}|d_E|\leq\gamma\sqrt{h}
\]
is complemented by two-sided volume density estimates and horizontal perimeter bounds in Carnot--Carath\'eodory balls centred on the interface, at radii up to a fixed multiple of $\sqrt{h}$. Their respective orders, $r^4$ and $r^3$, reflect the homogeneous dimension of $\mathbb{H}$. Combining these bounds with a volume-distance comparison and the discrete perimeter dissipation yields compactness as $h\to0$.

Our main result, Theorem~\ref{thmflatflow}, gives a sequence $h_i\to0$ along which the discrete evolutions converge in $L^1$ at every positive time to a bounded set $E(t)$ of finite horizontal perimeter. The limiting evolution satisfies
\[
|E(t)\triangle E(t')|\leq C|t-t'|^{1/2},
\qquad t,t'>0,\quad |t-t'|<1,
\]
where $C$ depends only on the initial horizontal perimeter and the isoperimetric constant. If the initial set also satisfies $|\partial E|=0$, the same estimate holds up to $t=0$, with $E(0)=E$. We call such a subsequential limit a flat flow. Its uniqueness and its agreement with smooth horizontal mean curvature flow or with the viscosity level-set evolution are not established here. 

The paper is organized as follows. Section~2 recalls the geometry of $\mathbb{H}$ and its distance functions. Section~\ref{sec3} collects the relevant facts about horizontal perimeter and mean curvature and develops the calibration argument for outward minimizing sets. Section~4 introduces the scheme, proves the density and displacement estimates, and constructs the flat flow, before discussing the remaining consistency questions.

{\bf{Aknowledgments}} A.C and M.N. are members of the Gruppo Nazionale per l'Analisi Matematica,
la Probabilit\`a e le loro Applicazioni (GNAMPA) of the Istituto Nazionale di Alta Matematica (INdAM).

\section{The Heisenberg group}
 
We briefly recall relevant definitions  and refer to  \cite[Chapter 2]{CapognaDanielli} for more details.
The Heisenberg group $\mathbb{H}$ is the Lie group $\R^3$ with the usual differentiable structure and composition given by 
\[(x_1,y_1,z_1) \opG (x_2,y_2,z_2)= (x_1+x_2,y_1+y_2,z_1+z_2+\frac{1}{2} (x_1y_2- x_2y_1)).\]  
The dilations of parameter $\lambda>0$ in $\mathbb{H}$ are defined as 
\[\delta_\lambda(p)=\delta_\lambda(x,y,z)=(\lambda x, \lambda y, \lambda^2z). \] The Haar measure on $\mathbb{H}$ is the Lebesgue measure in $\R^3$.
The left translation $L_p:\mathbb{H} \to \mathbb{H}$ is defined simply as $L_p(q)=p\opG q=(x_1+x_2,y_1+y_2,z_1+z_2+\frac{1}{2} (x_1y_2- x_2y_1))$, where $p=(x_1,y_1,z_1)$, $q=(x_2,y_2,z_2)$. 
 A straightforward calculation gives that for any $p$, $q$, $\bar{q}$ in $\mathbb{H}$, and $\lambda \in [0,1]$, we have that $L_{p}(\lambda q + (1-\lambda) \bar{q})= \lambda L_{p}(q) +(1-\lambda)L_{p}(\bar{q})$, yielding that Euclidean convex sets remain Euclidean convex after left-translation.
 
If $|E|$ denotes the Lebesgue measure of any measurable set $E$, we have  $|L_{p}(E)|=|E|$.  Moreover it holds $|\delta_{\lambda}(E)|=\lambda^{4}|E|$.
In the following we will denote
\[E^\lambda =\delta_\lambda ( E).\]

We introduce the  left invariant vector fields at every point $p=(x,y,z)\in \mathbb{H}$ 
 \[X(p)= \left(1, 0,  -\frac{y}{2}\right)^t \qquad Y(p)=\left(0,1,  \frac{x}{2}\right)^t\qquad T(p)= \left(0, 0,  1\right)^t\]
which will be identified with the differential operators 
\begin{align}\label{vectorfields} X(p)= \partial_x -\frac{y}{2}\partial_z\qquad Y(p)= \partial_y +\frac{x}{2}\partial_z\qquad T(p)= \partial_z.
\end{align}
The horizontal distribution $\mathcal{H}$  in every point $p=(x,y,z)\in \mathbb{H}$ is generated by the vector fields $X(p),Y(p)$, that is $\mathcal{H}(p)= $ Ker$ [\partial_z-\frac12 (x\partial_y-y\partial_x)]$. The horizontal distribution is nonintegrable since the Lie brackets of $X,Y$ satisfy
$[X(p), Y(p)]=T(p)$.  

A curve $\gamma:[0,1]\to \mathbb{H}$ is horizontal if at every $s \in [0,1]$, it holds
\[\gamma'_3(s)= -\frac{1}{2}\gamma_2(s)\gamma_1'(s)+\frac{1}{2}\gamma_1(s)\gamma_2'(s).  \]
Hence, given a planar curve $(\gamma_1, \gamma_2)$  in $\R^2$ and an initial point $z_0=\gamma_3(0)$,  there exists a unique lift  to a horizontal curve
  $(\gamma_1, \gamma_2, 
\gamma_{3}):[0,1]\to \mathbb{H}$ where
\begin{equation}\label{3componentecurva} \gamma_3(s)=\gamma_3(0)+\frac{1}{2} \int_0^s   -\gamma_2(s)\gamma_1'(s)+ \gamma_1(s)\gamma_2'(s) ds.\end{equation} 
In particular note that the height $|\gamma_3(1)-\gamma_3(0)|$ equals the  area   spanned by the planar  curve $(\gamma_1(s), \gamma_2(s))$ for $s\in [0,1]$.

We define the Carnot--Carathéodory distance among 2 points $p,q\in \mathbb{H}$ as follows:
\begin{multline}\label{dist1} d_{CC}(p,q)=\inf\{T\ \,|\, \exists \, \gamma:[0,T]\to \R^3, \text{s.t.} \gamma(0)=p, \gamma(T)=q,\\ (\gamma_1')^2+(\gamma_2')^2\leq 1\ \ \gamma_3'= \frac{1}{2}(-\gamma_2\gamma_1' +  \gamma_1\gamma_2')\}.\end{multline} 
It holds that $d_{CC}(p,q)=d_{CC}(w\opG p, w\opG q)$ for all $w\in \mathbb{H}$, so in particular $d_{CC}(p,q)=d_{CC}(q^{-1}\opG p, 0)=d_{CC}(0, p^{-1}\opG q)$ where $q^{-1}=(-x,-y, -z)$ if $q=(x,y,z)$, see \cite[Section 2.2]{CapognaDanielli}.

A sub-Riemannian metric on $\mathbb{H}$  is determined by any choice of inner product on
the horizontal subbundle $\mathcal{H}$. In this way it is possible to define  the length of horizontal curves and equip $\mathbb{H}$  with the structure of a metric
length space.  

We fix a  Riemannian metric $g$  on $\mathbb{H}$ such that  
the vector fields $X(p),Y(p),T(p)$ are orthonormal at every point $p=(x,y,z)$: 
\begin{equation}\label{matrix} g_p=\left(\begin{array}{ccc} 1+y^2/4 & -xy/4& y/2\\ -xy/4& 1+x^2/4& -x/2\\ y/2& -x/2& 1\end{array}\right).\end{equation} 
Note that $det g_p=1$, so the volume form coincides with the standard one.  
Note also that $X^t g_p= (1,0,0)$ and $Y^t g_p=(0,1,0)$, whereas $T^t g_p= (\frac{y}{2}, -\frac{x}{2}, 1)$.

We define the length of a horizontal curve $\gamma:[0,1]\to \mathbb{H}$ as
\[L(\gamma)=\int_0^1 \sqrt{(\gamma'(s)\cdot_g X(\gamma(s)))^2+ (\gamma'(s)\cdot_g Y(\gamma(s)))^2}ds= \int_0^1 \sqrt{ (\gamma'_1)^2+ (\gamma_2')^2}ds \]
where $\gamma'\cdot_g X(\gamma)= (\gamma')^t g_{\gamma} X(\gamma)$. 

\begin{remark}[Geodesics]\label{geodesic} 
\upshape  Recalling \eqref{3componentecurva} if  we fix $(0,0, 0)$ and $p=(x_1,y_1, z_1)$ in $\mathbb{H}$, the problem of finding a horizontal curve of minimal length connecting $0$ and $p$ coincides with the problem of  finding a curve in the plane of minimal length connecting $(0,0)$ and $(x_1,y_1)$ among the curves such that the signed area   of the region $S$ delimited by the curve and the segment joining $(0,0)$ and $(x_1,y_1)$ is fixed and equal to $z_1$.  Therefore a length minimizing curve between $0$ and $p$ is the lift of a
circular arc joining the origin  with $(x_1,y_1)$, whose convex hull has signed  area $z_1$. 
If $p=(0,0, z_1)$ there are infinitely many geodesics (which are all the lifts of circles of area $|z_1|$ passing through the origin), whereas if $x_1^2+y_1^2\neq 0$ it is possible to show that there exists a unique geodesic connecting $0$ and $p$ (see \cite[\S 2.3]{CapognaDanielli} for further details and a representation formula). \end{remark} 

\begin{definition}[Pansu sphere]\label{defpansu}
We consider the   geodesic $\gamma:[0,2\pi]\to \mathbb{H}$ joining $0$ and $(0,0, \pi R^2)$ which are parametrized as 
\[\gamma(s)=(R(1-\cos s), -R\sin s, R^2(s-\sin s)/2). \] The curvature of the planar  projection of the  geodesic is $1/R$. 
Rotating such a geodesic about the $z$-axis produces a surface of revolution
 which is called the Pansu sphere.
 \end{definition} The Pansu sphere is conjectured to be the isoperimetric set, see \cite{p}.

    \subsection{Metric structure} 
    \label{sec:distances}
In the previous section we introduced the Carnot--Carathéodory distance on $\mathbb{H}$. There are  other possible definitions of distance. 
For $p,q \in \mathbb{H}$ we have: 
\begin{itemize}
\item Carnot--Carathéodory distance \begin{equation}\label{dcc} d_{CC}(p,q) = \inf \{ L(\gamma) \ |  \ \gamma:[0,1] \to \mathbb{H}^{1} \text{ is horizontal}, \gamma(0)=p, \gamma(1)=q \};\end{equation}   This distance is   left invariant and $1$-homogeneous with respect to dilations. 
\item $d_{\infty}(p,q):= \| p^{-1} \opG q\|_{\infty}$, where the box norm $\|\cdot \|_{\infty}$ is given by  $$\|p \|_{\infty}= \max\{(x^2+y^2)^{\frac12}, 2|z|^{\frac{1}{2}} \};$$
\item $d_{Kor}(p,q):= \| p^{-1} \opG q\|_{Kor}$, where the Kor\'{a}nyi norm $\| \cdot\|_{Kor}$ is given by  $$\|p \|_{Kor}= (|x^2+y^2|^{2}+ 16 |z|^2 )^{1/4}.$$
\end{itemize}
The above distances are all equivalent (\cite[page~5]{Monti}),   i.e. there exists an absolute constant $C>0$ such that
\begin{align}\label{const-disteq}
C d_{\cdot}(p,q) \leq d_{CC}(p,q) \leq C^{-1}d_{\cdot}(p,q) \qquad \forall p,q \in \mathbb{H},  \qquad  \text{ for } \quad d_{\cdot}=d_{\infty} \text{ or } d_{\cdot}=d_{Kor},
\end{align}
 hence they induce the same topological structure. 

Based on this equivalence, we have the following relation that will be useful later on: for every $\lambda>0$, 
\begin{equation}\label{dcclambda} 
d_{CC}(p, \delta_\lambda p)\leq C\|p\|_{\infty}  \max\left(|\lambda -1|, \sqrt{|\lambda^{ 2}-1|} \right). \end{equation} 
Indeed by $1$-homogeneity and \eqref{const-disteq} we have, renaming the constant $C$ from line to line, 
\begin{align*}  d_{CC}(p, \delta_{ \lambda} p)& =d_{CC} (0, p^{-1}\opG \delta_{\lambda} p)  \leq  C\|p^{-1} \opG \delta_{\lambda} p\|_\infty=
C\|(\lambda -1)x, (\lambda -1)y, (\lambda^{2}-1)z)\|_\infty\\& =C\max\left(|\lambda -1|\sqrt{x^2+y^2}, 2\sqrt{|\lambda^{2}-1|}\sqrt{|z|}\right)\leq C\|p\|_{\infty}  \max\left(|\lambda -1|, \sqrt{|\lambda^{ 2}-1|} \right)\end{align*} where $p=(x,y,z)$.

  The relation to the Euclidean one is explained in the following
\begin{proposition}{\cite[Prop.2.1]{serra}}\label{serra21}
For any bounded subset $\Omega \subset \mathbb{H}$ there exist positive constants $c_{1}(\Omega)$, $c_{2}(\Omega)$ such that
$$ c_{1}(\Omega) |p-q|_{\R^{3}} \leq d_{CC}(p,q) \leq c_{2}(\Omega) |p-q|^{\frac{1}{2}}_{\R^{3}}$$
for $p,q \in \Omega$. In particular, the topologies defined by $d_{CC}$ and by the Euclidean
distance coincide.
\end{proposition}
It is possible to associate to every distance above the corresponding family of balls. In  particular we will consider the $CC$ balls 
\begin{equation}\label{ccball}
B_{CC}(p,r)=\{  q \in \mathbb{H}\,:\, d_{CC}(p,q) <r\}.
\end{equation}
Note that $$ |B_{CC}(p,r)|=r^{4}|B_{CC}(0,1)|$$
for any $p\in \mathbb{H}$ and $r>0$.
For $p \in \mathbb{H}$ and  sets  $A, B \subset \mathbb{H}$ let us introduce the distance between $p$ and the set $A$ and the Hausdorff distance between $A,B$:  \begin{equation}\label{distanzeinsiemi} \text{dist}(p,A):=\inf_{q\in  A} d_{CC}(p,q)\qquad \text{dist}(A,B)=\max\left\{\sup_{a\in A} \text{dist}(a,B), \sup_{b\in B} \text{dist}(b,A)\right\}.\end{equation}
  
 Given $E \subset \mathbb{H}$ open bounded set we define the \emph{signed distance}  to be    
\begin{equation} \label{sig}
d_E(p)= \begin{cases}d_{CC}(p,  E), &p \in \mathbb{H} \setminus E, \\
 - d_{CC}(p, \mathbb{H}\setminus   E), & p \in E. \end{cases}
\end{equation}
The signed distance function is a Lipschitz map with respect to the metric structure in $\mathbb H$, $|d_{E}(p)-d_{E}(q)| \leq C d_{CC}(p,q)$, hence only locally H\"older continuous with respect to the Euclidean distance by Prop. \ref{serra21} and the equivalence of the above norms \cite{Monti}.

 The metric space $(\mathbb{H}, d_{CC})$  is complete and locally compact. By the definition of $d_{CC}$, it
is also a length space, and also  a geodesic space, namely for all $p, q\in\mathbb{H}$ there exists a horizontal curve $\gamma$ with $\gamma(0)=p$ and $\gamma(1)=q$, such that $d_{CC}(p,q)=L(\gamma)$, see \cite{Monti}, \cite[Chapter 2]{CapognaDanielli}. 

 For any measurable set $F$ and $t \in [0,1]$ we define $F(t)$ to be the set of points where the volume density of $F$ equals $t$, i.e.
$$ F(t)= \left \{ p \in \mathbb{H} \ :\, \lim_{r\to 0} \frac{|F \cap B_{CC}(p,r)|}{|B_{CC}(p,r)|}=t \right \}$$
In the following we  identify any measurable set with the representative given by the set of Lebesgue points of the characteristic functions, i.e. with $F(1)$.

 The essential boundary $\partial^{e}F$ of $F$ is defined by $\partial^{e}F:= \mathbb{H} \setminus (F(1) \cup F(0))$.  
 A set $B$ is said to be $C^{k}$ if its boundary $\partial B$ is $C^{k}$ in the classical (euclidean) sense. 
 
We introduce the definition of sets which are strictly starshaped with respect to the origin, which will be useful to construct foliations of the space $\mathbb{H}$, see \cite{fog,dragoni}. 

\begin{definition}\label{strictlystarshaped} We say that a $C^{1}$
   bounded open set $B$ which contains $0$,  is strictly starshaped with respect to $0$ (with respect to the group of dilations)  if \begin{equation}\label{starshaped}\forall p\in \partial B\quad Z(p)\cdot N(p)>0,\end{equation} where $N(p)$ is the outer Euclidean normal, and for $p=(x,y,z)$,  $Z(p)=(x,y,2z)$ is the vectorial field generated by the dilations. \end{definition} 
 In particular strictly starshaped sets fibrate all the space by dilations $B^\lambda=\delta_\lambda B$. 
\begin{lemma}\label{lemmastellati} 
Let $B$ a $C^{1}$ bounded open  set  which contains $0$ and which is strictly starshaped with respect to the origin according to the previous definition. 
Then $\mathbb H=\cup_{\lambda\geq 0}B^\lambda$ and  for every $p\in\mathbb{H}\setminus\{0\}$ there exists a unique $\lambda$ such that $p\in \partial B^\lambda$. 
 \end{lemma} 
 \begin{proof}
 Obviously $\delta_0 p=0$ for every $p$. 
 
 Let us assume that $B=\{p, \ F(p)<0\}$, $\partial B=\{p, \ F(p)=0\}$  for some $C^{1}$   function $F:\mathbb H\to \R$ with $|\nabla F|\neq 0$ on $\partial B$. Therefore the condition that $B$ is strictly starshaped reads  $\nabla F(p) \cdot Z(p)>0$ for all $p\in \partial B$.
 
First of all we observe  that for every $q\in \mathbb H\setminus \{0\}$ there exist  $\lambda>0$ and $p\in \partial B$ such that $\delta_\lambda p=q$. Let us define $h(\lambda):= \delta_\lambda q$ for $\lambda\in [0, +\infty)$. Then $h(0)=0\in B$ and $h(1)=q$. Moreover since $B$ is bounded there exists $\lambda_0>0$ such that $h(\lambda_0)\not\in B$. 
 Therefore $F(h(0))<0$, whereas $F(h(\lambda_0))>0$, and the function $\lambda\to F(h(\lambda))$ is continuous: this implies that there exists $\bar \lambda\in (0,\lambda_0)$ such that $F(h(\bar \lambda))=0$. Hence  $\delta_{\bar \lambda } q=p\in \partial B$, and $q=\delta_{\frac{1}{\bar\lambda}}p$. 

Next we prove that such $\lambda>0$ and $p\in \partial B$ are unique. Assume by contradiction that there exists $p_1, p_2\in \partial B$, $\lambda_1, \lambda_2$ such that $q=\delta_{\lambda_1}p_1=\delta_{\lambda_2}p_2$. Without loss of generality we may fix  $0<\frac{\lambda_1}{\lambda_2}\leq 1$. Note that $\delta_{\frac{\lambda_1}{\lambda_2}}p_1=p_2$. 
 
 Now we claim  that for every $p\in \partial B$ it holds that  $\delta_\lambda p\in B$ for all $\lambda<1$, that is $F(\delta_\lambda p)<0$.  If the claim is true, we conclude from $\delta_{\frac{\lambda_1}{\lambda_2}}p_1=p_2$ with  $p_1, p_2\in \partial B$  and  $\frac{\lambda_1}{\lambda_2}\leq 1$ that    necessarily $\frac{\lambda_1}{\lambda_2}=1$, and then also $p_1=p_2$.

We are left to prove the claim. Let us define $g(\lambda):= F(\delta_\lambda p)$. 
We have that $g(1)=0$ and $g(0)<0$ (since $0$ is in  $B$). We compute  $g'(\lambda)= \nabla F(\delta_\lambda p) \cdot \frac{1}{\lambda} Z(\delta_\lambda p)$. In particular $g'(1)>0$. Therefore there exists $\hat{\delta}>0$ such that $g(\lambda)<0$ for every $\lambda \in (1-\hat{\delta}, 1)$. Let's assume by contradiction that there exists  $\lambda_0=\sup\{\lambda\leq 1\ \text{ s.t. } g(\lambda)\geq 0\}$. Then $0<\lambda_0\leq 1-\hat{\delta}$. By continuity $g(\lambda_0)=0$ and then $\delta_{\lambda_0}p\in\partial B$, and $g(\lambda)<0$ for $\lambda\in (\lambda_0, 1)$. Recall that  
$g'(\lambda_0)=  \nabla F(\delta_{\lambda_0} p) \cdot \frac{1}{\lambda_0} Z(\delta_{\lambda_0}p)>0$ 
 since $\delta_{\lambda_0} p\in \partial B$. This contradicts the fact that $g(\lambda_0)=0$ and $g(\lambda)<0$ for $\lambda>\lambda_0$ and so, we conclude  that  
 for every $p\in \partial B$, it holds that $\delta_\lambda p\in B$ for every $\lambda\in [0,1)$. 
 
\end{proof} 
  \subsection{Differential structure} 

We will indicate the Levi-Civita connection in $(\mathbb{H},g)$ as $D$, see \cite[Chapter 2]{CapognaDanielli}. Upon recalling \eqref{vectorfields} we have that \begin{eqnarray}   D_X(X)=0 &D_X (Y)=\frac{1}{2}T & D_X(T)=-\frac{1}{2}Y\\ D_Y(X)=-\frac{1}{2}T &D_Y (Y)=0 & D_Y(T)=\frac{1}{2}X\\  D_T(X)=-\frac{1}{2}Y &D_T (Y)=\frac{1}{2}X & D_T(T)=0. 
\end{eqnarray} 
For a vector field $U$ we define \begin{equation}\label{J} J(U):=D_U (2T).\end{equation}  Hence by the Lie bracket condition we have that $J(X)=-Y$, $J(Y)=X$  and $J(T)=0$.  
Thus $J$ defines an isometry on the horizontal distribution and $J^2=-Id$. 

Given a $C^1$ function $f:\mathbb{H}\to \R$ we define the horizontal gradient at $p=(x,y,z)$ as
\begin{equation}\label{horgrad}
\nabla_{\mathbb{H}} f =Xf X+Yf Y=  (\partial_xf-\frac{y}{2}\partial_z f)X+(\partial_y f+\frac{x}{2}\partial_z f)Y.
\end{equation}
We point out that 
\[\nabla_{\mathbb{H}} f (p)= \left(\begin{array}{ccc}1 &0 & -y/2\\ 
0& 1& x/2\end{array}\right)\cdot \nabla f(p)\]
and that we are identifying a vector $\xi=(\xi_1, \xi_2)$ in the horizontal distribution $\mathcal{H}(p)$ with a vector $\tilde \xi\in\R^3$ through
\[\tilde \xi=\left(\begin{array}{cc}1 &0 \\ 0 &1\\  -y/2 & x/2\end{array}\right)\cdot \xi. \]

Given a $C^1$ vector field $U:\mathbb{H}\to \mathbb{H}$, we may express it   using both the
basis $(X,Y,T)$ and the standard basis of vector fields   in $\R^3$ (with $e_3=T$) at every point $p\in \mathbb{H}$: 
\[U(p)= u_1X+u_2Y+u_3T = u_1(1,0,0)^t +u_2(0,1,0)^t+ \left(-\frac{y}{2}u_1+\frac{x}{2}u_2+ u_3\right)(0,0,1)^t.  \]

We define the    horizontal projection of $U$  on the horizontal distribution $\mathcal{H}$, that is: at $p=(x,y,z)$, we have
\begin{equation} U_{\mathbb{H}}(p) = u_1X+u_2Y. \label{horizontalprojection} \end{equation}

\begin{definition} We have that a vector field is horizontal iff $U=U_\mathbb{H}$, that is $  u_3=0$.
In particular observe that every   function in  $C^1(\R^3, \R^2)$ can be associated in a unique way to a horizontal vector field in $\mathbb{H}$. 
\end{definition}

 We define the  divergence in $\mathbb{H}$ as follows: given a $C^1$ vector field we have that 
\[\Div_g  U= D_X U\cdot_g X+ D_Y U\cdot_g Y+ D_T U\cdot_g T.\] Note that $  D_XU\cdot_g X= D_XU\cdot_g X+U\cdot_g D_X X= X(U\cdot_g X)= X(u_1)$, and analogously 
$  D_YU\cdot_g Y= Y(u_2)$ and $  D_TU\cdot_g T=T(u_3)$. Therefore 
\[\Div_g  U= X(u_1)+Y(u_2)+T(u_3)=
\partial_x u_1+\partial_y u_2 -\frac{y}{2}\partial_z u_1+\frac{x}{2}\partial_z u_2+\partial_z u_3.\]
 The horizontal divergence of a vector field will be given by 
 \begin{equation}\label{hordiv}\Div_{\mathbb{H}} U= X(u_1)+Y(u_2)\end{equation} and it coincides with the divergence of its horizontal projection
\[ \Div_{\mathbb{H}} U=\Div_g  U_{\mathbb{H}}=\partial_x u_1+\partial_y u_2 -\frac{y}{2}\partial_z u_1+\frac{x}{2}\partial_z u_2.\]
In particular if a vector field is horizontal we have that $\Div_g U=\Div_{\mathbb H}U$. 
In terms of the standard Euclidean divergence $\Div_{\R^{3}}$ observe that
\[
\Div_{\R^{3}} U_{\mathbb{H}} =\Div_{\R^{3}} (u_{1}X + u_{2}Y) = X(u_{1})+Y(u_{2})
\]
therefore we can also state that
\begin{align}\label{div-equi}
 \Div_{\R^{3}} U_{\mathbb{H}}=\Div_{\mathbb{H}} U_{\mathbb{H}}=\Div_{g} U_{\mathbb{H}}.
\end{align}
In the following we denote the Euclidean divergence  $\Div_{\R^{3}}$ simply by $\Div$.

 \section{Horizontal perimeter and  mean curvature} \label{sec3}
We collect here the most relevant  notions on the geometry of manifolds in $\mathbb{H}$; for more details we refer to  \cite[Chapter 4]{CapognaDanielli} and \cite{rr}. 
We consider a  $C^1$ surface $\Sigma$  enclosing an open set $E$, i.e.
 $$\partial E=\Sigma.$$
 Moreover, we  assume that it has regular points, that is, it can be locally written as the $0$ level set of a $C^1$ function with non vanishing gradient, and   $E=\{p\ |\ f(p)<0\}$ for some $f\in C^1$.
 
We denote by $N(p)=\frac{\nabla f(p)}{|\nabla f(p)|}$ the euclidean exterior  normal to the surface $\Sigma$ at a point $p=(x,y,z)$. We recall that at $p$ there holds
 \[ Xf= \partial_x f-\frac{y}{2} \partial _z f \qquad Yf= \partial_y f+\frac{x}{2} \partial _z f  \qquad Tf=   \partial _z f.\]  Following  \cite[Chapter 4]{CapognaDanielli}, we define: 
\begin{eqnarray}\label{e1}e_1:&=& \frac{Yf}{\sqrt{(Xf)^2+(Yf)^2}} X- \frac{Xf}{\sqrt{(Xf)^2+(Yf)^2}} Y\\\nonumber  
e_2:&=& \frac{Tf}{\sqrt{(Xf)^2+(Yf)^2+(Tf)^2}}\frac{Xf}{\sqrt{(Xf)^2+(Yf)^2}} X+ \frac{Tf}{\sqrt{(Xf)^2+(Yf)^2+(Tf)^2}} \frac{Yf}{\sqrt{(Xf)^2+(Yf)^2} } Y 
\\ & & -  \frac{\sqrt{(Xf)^2+(Yf)^2}}{\sqrt{(Xf)^2+(Yf)^2+(Tf)^2}} T.\label{e2}
\end{eqnarray} 
We observe that  \[ Xf = X(p)\cdot \nabla f(p)\qquad Yf=Y(p)\cdot \nabla f(p) \qquad Tf= T(p)\cdot \nabla f(p).\]  Using this it is 
 straightforward to check that
\begin{eqnarray*}  e_1\cdot N(p)&=& \frac{1}{  \sqrt{(Xf)^2+(Yf)^2}|\nabla f|} \left[Yf (X\cdot \nabla f)- Xf (Y\cdot \nabla f)\right]=0\\ 
 e_2\cdot N(p)&=&  \frac{1}{  \sqrt{(Xf)^2+(Yf)^2 +(Tf)^2}|\nabla f|} \left[Tf \frac{(Xf)^2}{ \sqrt{(Xf)^2+(Yf)^2}}+ Tf \frac{(Yf)^2}{{ \sqrt{(Xf)^2+(Yf)^2}}} + \right.\\ &&\left. -\sqrt{(Xf)^2+(Yf)^2} Tf\right]=0. 
\end{eqnarray*} 
Therefore $\{e_1, e_2\}$ is a basis of  the tangent space to $\Sigma$ at the point $p$. Moreover one can verify that  the basis  $\{e_1, e_2\}$ is orthonormal with respect to the Riemannian metric $g$, that is 
$e_1 \cdot_g e_1=1=e_2 \cdot_g e_2$, and $e_1 \cdot_g e_2=0$. 

We define the Riemannian normal  at $p$ to $\Sigma$ as follows:
\[N_g(p)=\frac{Xf}{\sqrt{(Xf)^2+(Yf)^2+(Tf)^2}} X+ \frac{Yf}{\sqrt{(Xf)^2+(Yf)^2+(Tf)^2}}  Y +\frac{Tf}{\sqrt{(Xf)^2+(Yf)^2+(Tf)^2}} T.\] 
Note that $N_g\cdot_g e_1=0= N_g\cdot_g e_2$.  So $(N_g, e_1, e_2)$ form a basis of the tangent bundle of $(\mathbb{H}, g)$ at $p$, whereas $(e_1, e_2)$ form a basis of the tangent bundle of $\Sigma$ at $p$.
Let 
\begin{equation}\label{nh}N_{\mathbb{H}}(p)=N_g(p)- (N_g(p)\cdot_g T(p))T(p) = \frac{Xf}{\sqrt{(Xf)^2+(Yf)^2+(Tf)^2}} X+ \frac{Yf}{\sqrt{(Xf)^2+(Yf)^2+(Tf)^2}}  Y .  \end{equation} 

\begin{definition}The singular points in $\Sigma$ are the points $p$ such that $N_{\mathbb{H}}(p)= 0$, that is the points where $N(p)\cdot X(p)=0$ and $N(p)\cdot Y(p)=0$.  This implies that 
$N(p)$ is orthogonal to the horizontal bundle $\mathcal{H}(p)$ and therefore the horizontal bundle $\mathcal{H}(p)$ and   the euclidean tangent space to $\Sigma$ at $p$ coincide.  
We denote  the set of  singular points  as $\Sigma_0$.   \end{definition} 
\begin{remark}\upshape  The set $\Sigma_0$ is closed and has empty interior in  $\Sigma$.  Moreover, for a $C^2$ surface  $\Sigma$, it is known that the Hausdorff dimension of  $\Sigma_0$ with respect to the Riemannian distance in  $\mathbb H$ 
is less than or equal to $1$.  We refer to \cite[Section 2]{rr} for more details and references.
\end{remark}

The \emph{exterior horizontal normal} to $\Sigma$ at $p\in \Sigma\setminus \Sigma_0$ is the  normalization of $N_{\mathbb {H}}$ with respect to the $g$ metric: 
\begin{equation}\label{normal}
\nu_{\mathbb{H}}(p)=  \frac{Xf}{\sqrt{(Xf)^2+(Yf)^2}} X+ \frac{Yf}{\sqrt{(Xf)^2+(Yf)^2}}  Y = \frac{\nabla_{\mathbb H}f(p)}{|\nabla_{\mathbb H} f(p)|_g}=\frac{ N_{\mathbb H}(p)}{|N_{\mathbb H}|_g(p)}.
\end{equation}

Note that for every $p\in \Sigma\setminus \Sigma_0$   \begin{eqnarray}\label{e1nuovo} e_1(p)&=& J(\nu_{\mathbb{H}}(p)) \\ \label{e2nuovo} e_2(p)&=& (N_g\cdot_g   T)\  \nu_{\mathbb H}(p)- |N_{\mathbb H}|_g T .\end{eqnarray} In particular $\{e_1, \nu_{\mathbb{H}}\}$ forms an orthonormal basis of the horizontal bundle $\mathcal{H}(p)$ at every $p\in \Sigma\setminus \Sigma_0$. 
 
 \begin{remark}[Horizontal normal of planes]\upshape\label{planenormal}
We consider the   euclidean plane $\{(x,y,z)\ | \ z=k\}$. It is easy to see that the only characteristic point is $(0,0, k)$ and that the horizontal normal  coincides with \[ \nu_{\mathbb{H}}(p)= \frac{-y}{\sqrt{y^2+x^2}}X+ \frac{x}{\sqrt{y^2+x^2}}Y. \]
More generally if we consider a non-vertical  plane   $\{(x,y,z)\ | \ a x+by+z=k\}$ then the only characteristic point is $(-2b, 2a, k)$ and the horizontal normal is given by
\[ \nu_{\mathbb{H}}(p)= \frac{(a-y/2)X+ (b+x/2)Y}{\sqrt{(a-y/2)^2+(b+x/2)^2}}. \]
Finally   the vertical planes  $\{(x,y,z)\ | \ a x+by=k\}$ have no singular points and  \[\nu_{\mathbb{H}}(p)= \frac{aX+bY}{\sqrt{a^2+b^2}}. \]
\end{remark}

\subsection{Horizontal mean curvature} 
To introduce the notion of mean curvature we assume now that the surface $\Sigma=\partial E$ is of class $C^2$. 
 Recalling \eqref{e1nuovo}, we have that  
 \[\{\nu_{\mathbb{H}}(p), e_1(p), T\}\]  forms an orthonormal basis of the tangent bundle of the Riemannian space $(\mathbb{H}, g)$ at every  $p\in \Sigma\setminus \Sigma_0$.

Moreover, we already observed that $e_1(p)$ belongs to the  tangent space to $\Sigma$. The integral curves of $e_1$ in $\Sigma\setminus \Sigma_0$ are  
called characteristic curves and foliate the regular part of $\Sigma$. 
 
 Since $|\nu_{\mathbb{H}}|_g=1$,  it holds $D_{\nu_{\mathbb{H}}}\nu_{\mathbb{H}}\cdot_g \nu_{\mathbb{H}}=0$, and since  $\nu_{\mathbb{H}}\cdot_gT=0$ and  $D_TT=0$, we get that
$D_T \nu_{\mathbb{H}}\cdot_g T=- D_T T\cdot_g \nu_{\mathbb{H}}=0$. Therefore 
\begin{equation}\label{divnormal}\Div_g \nu_{\mathbb{H}}(p)=D_{ \nu_{\mathbb{H}}}\nu_{\mathbb{H}} \cdot_g \nu_{\mathbb{H}}+D_{T}\nu_{\mathbb{H}} \cdot_g T+ D_{e_1}\nu_{\mathbb{H}} \cdot_g e_1= D_{e_1}\nu_{\mathbb{H}} \cdot_g  e_1= -D_{e_1}e_1\cdot_g\nu_{\mathbb{H}} .   \end{equation}
Recall that an equivalent formulation of this divergence using the basis $\{X,Y,T\}$ is given by
\begin{eqnarray}\label{eq21} \nonumber \Div_g \nu_{\mathbb{H}}(p) &=&  D_{X}\nu_{\mathbb{H}} \cdot_g X+D_{T}\nu_{\mathbb{H}} \cdot_g T+ D_{Y}\nu_{\mathbb{H}} \cdot_g Y\\ \label{divnormal2}&=& X\left(\frac{Xf}{\sqrt{(Xf)^2+(Yf)^2}}\right)+Y\left(\frac{Yf}{\sqrt{(Xf)^2+(Yf)^2}} \right)=\Div_\mathbb{H} \nu_{\mathbb{H}}(p).   \end{eqnarray}
 We consider the  vector $e_2$ in the tangent bundle,  defined in \eqref{e2}. As observed in the previous section   the tangent bundle of $\Sigma$ at $p$ has orthonormal basis given by $\{e_1, e_2\}$.
Then the Riemannian mean curvature of $\Sigma$ at $p \in \Sigma $ is defined as the tangential divergence of the normal $N_g$:
\begin{equation}\label{tangentialdivergence}H_g(p, \Sigma)=\Div_{\Sigma}(N_g)= D_{e_1} N_g\cdot_g e_1+ D_{e_2} N_g\cdot_g e_2.\end{equation} 

We compute now the tangential divergence of the horizontal normal $\nu_{\mathbb H}$, Note that $D_{e_1}T=\frac{1}{2} J(e_1)=-\frac12 \nu_{\mathbb H}$, therefore $D_{e_1}T \cdot_g  e_1=0$.  
So, recalling the definition of $\nu_{\mathbb H}$, we get 
 \[  |N_\mathbb H|_g D_{e_1}\nu_{\mathbb{H}} \cdot_g  e_1=  D_{e_1} N_g \cdot_g  e_1-(T\cdot_g N_g) (D_{e_1}T \cdot_g  e_1)=  D_{e_1} N_g \cdot_g  e_1. \]
On the other hand, since $\{e_1, T, \nu_{\mathbb H}\}$ forms an orthonormal basis of the tangent bundle of   $(\mathbb{H}, g)$, and $D_{e_2}\nu_{\mathbb H}\cdot_g \nu_{\mathbb H}=0$ (due to the fact that $|\nu_{\mathbb H}|_g=1$) we have that 
\[ D_{e_2} \nu_{\mathbb H}=( D_{e_2} \nu_{\mathbb H}\cdot_g e_1)e_1+ (D_{e_2} \nu_{\mathbb H}\cdot_g T)T. \]
Observe that, due to the definition of $J$, 
\[D_{e_2} \nu_{\mathbb H}\cdot_g T= - D_{e_2} T \cdot_g  \nu_{\mathbb H}= -\frac{1}{2} J(e_2)\cdot_g  \nu_{\mathbb H}= -\frac{1}{2} e_2\cdot_g  J(\nu_{\mathbb H})=-\frac12 e_2\cdot_g e_1=0. \] Therefore  $D_{e_2} \nu_{\mathbb H}=( D_{e_2} \nu_{\mathbb H}\cdot_g e_1)e_1$, and then
\[\Div_{\Sigma}(\nu_{\mathbb H})= D_{e_1} \nu_{\mathbb H}\cdot_g e_1+ D_{e_2} \nu_{\mathbb H}\cdot_g e_2=D_{e_1}\nu_{\mathbb{H}} \cdot_g  e_1=  \frac{1}{|N_{\mathbb H}|_g} D_{e_1} N_g \cdot_g  e_1.\]

We define the horizontal mean curvature of $\Sigma$ at $p \in \Sigma\setminus \Sigma_0 $ as the  tangential divergence of the horizontal  normal, written above, that is 
\begin{equation}\label{horizontalmc}H_{\mathbb{H}}(p, E )=\Div_{\Sigma}(\nu_{\mathbb H})=  \frac{1}{|N_{\mathbb H}|_g} D_{e_1} N_g \cdot_g  e_1=D_{e_1}\nu_{\mathbb{H}} \cdot_g  e_1= \Div_g  \nu_{\mathbb{H}}(p).\end{equation} It coincides with the divergence of $\nu_{\mathbb H}$ (see \eqref{eq21}). 

\begin{remark}\upshape\label{remarkhomoH} 
The horizontal mean curvature is $-1$ homogeneous with respect to dilations and constant with respect to left translations, i.e.
$H_{\mathbb{H}}(\delta_{\lambda}(p), \delta_{\lambda}(E))=\frac{1}{\lambda}H_{\mathbb{H}}(p, E ) $ and that for every $q\in \mathbb H$, $H_{\mathbb{H}}(L_{q} (p), L_{q}( E))= H_{\mathbb{H}}(p, E ) $. \end{remark} 
 We recall from \cite{af1},  \cite{af2}, the following result relating the horizontal gradient and horizontal divergence of the signed distance function from a smooth set to the horizontal normal and horizontal curvature. 
\begin{theorem}\label{teo:7.1} Let $E\subseteq \mathbb{H}$ be a bounded $C^k$ open set (in Euclidean sense), for $k\geq 2$,  and let $d_E$ be the signed distance function from $E$ as defined in \eqref{sig} (with respect to the Carnot--Carathéodory distance $d_{CC}$). Then the following hold:
\begin{enumerate}
\item $d_E$ and $\nabla_{\mathbb H}d_E$ are of class $C^{k-1}$ in Euclidean sense in an open neighborhood of $\partial E\setminus \Sigma_0$, where $\Sigma_0$ is the set of singular points in $\partial E$. Moreover at   every $p\in\partial E\setminus \Sigma_0$, there holds  $\nabla_{\mathbb H} d_E(p)=\nu_{\mathbb{H}}(p)$ and $|\nabla_{\mathbb H} d_E(p)|=1$. 
\item Let $k\geq 3$. At every $p\in\partial E\setminus \Sigma_0$, it holds that \[\nabla_{\mathbb H} d_E(p)=\nu_{\mathbb{H}}(p) \quad\text{ and}\quad \Div_{\mathbb H}\nabla_{\mathbb H} d_E(p)= H_{\mathbb H}(p, E).\] 
\end{enumerate}
\end{theorem}
\begin{proof} For the proof   we refer to \cite[Theorem 1.1, Corollary 6.2]{af1}, \cite[Theorem 1.1, Corollary 4.1]{af2}. \end{proof} 

\begin{remark}\label{HXLS}\upshape
Using \eqref{horizontalmc} and \eqref{eq21}, the horizontal mean curvature of $\Sigma=\{ f=0 \}$, boundary of $E=\{f\leq 0\}$,  at $p \in \Sigma\setminus \Sigma_0 $ can be expressed as
\begin{align*}
H_{\mathbb{H}}(p, E)=X\left(\frac{Xf}{\sqrt{(Xf)^2+(Yf)^2}}\right)+Y\left(\frac{Yf}{\sqrt{(Xf)^2+(Yf)^2}} \right)\\
=\frac{1}{|\nabla_{\mathbb{H}}f|_{g}} tr \left[\left(I_{2}-\frac{\nabla_{\mathbb{H}}f \otimes \nabla_{\mathbb{H}}f}{|\nabla_{\mathbb{H}}f|_{g}^{2}}\right) (\nabla^{2}_{\mathbb{H}}f)^{*}  \right],
\end{align*}
where $I_{2}$ is the identity matrix in $\R^{2}$, $|\nabla_{\mathbb{H}}f|_{g}=\sqrt{(Xf)^{2}+(Yf)^{2}}$, $\nabla_{\mathbb{H}}f=(Xf,Yf)$ (recall \eqref{horgrad}) and $(\nabla^{2}_{\mathbb{H}}f)^{*}$ is the matrix
\begin{align*}
(\nabla^{2}_{\mathbb{H}}f)^{*}:=\left( \begin{array}{cc}
XXf & \frac{1}{2} (XYf + YXf)\\
\frac{1}{2} (XYf + YXf) & YYf
 \end{array} \right).
\end{align*}
Moreover, with the notation  $X_{1}:=X$, $X_{2}:=Y$ we can also write
\begin{align*}
H_{\mathbb{H}}(p, E ) |\nabla_{\mathbb{H}}f|_{g}= \sum_{i,j=1}^{2}\left(\delta_{ij} -\frac{(X_{i}f) ( X_{j}f)}{|\nabla_{\mathbb{H}}f|^{2}_{g}} \right)X_{i}X_{j}f.
\end{align*}
\end{remark}

For axially symmetric sets we have a more  explicit formula for the horizontal mean curvature.
\begin{remark}[Horizontal mean curvature for axially symmetric sets] \upshape \label{remarkaxi} 
Let $\Sigma$ be  an axially symmetric surface, that is $\Sigma=\{(x,y,z) \ | \ |z|= f(\sqrt{x^2+y^2})\}$ for some $f\in C^2$. The singular set is given by $\Sigma_{0}\subseteq \{(0,0,\pm f(0)) \}$. The horizontal outer  normal to $E=\{(x,y,z) \ | \ |z|<f(\sqrt{x^2+y^2})\}$ can be easily computed at a point $p=(x,y,z)$, with $r=\sqrt{x^2+y^2}>0$ and $z>0$,  as 
\begin{equation}\label{radialnormal}\nu_{\mathbb{H}}(p)= \frac{1}{\sqrt{(f'(r))^2+r^2/4}}\left[\left(-f'(r)\frac{x}{r}-\frac{y}{2}\right)X+ \left(-f'(r)\frac{y}{r}+\frac{x}{2}\right)Y\right].\end{equation}

So the horizontal mean curvature at $p\in \Sigma\setminus \Sigma_0$ (with $z>0$) can be computed as follows (for details see \cite[Section 4.3]{CapognaDanielli}):   
\begin{equation}\label{radialmc}H_{\mathbb{H}}(p, E )=- \frac{\frac14r^2f''(r)+\frac{1}{r}(f'(r))^3}{[(f'(r))^2+r^2/4]^{3/2}}  \end{equation}
For a comparison, the euclidean counterparts, again at  a point $p=(x,y,z)$, with $r=\sqrt{x^2+y^2}>0$ and $z>0$ are given by the euclidean outer normal
\begin{align}\label{radnormalEuc}
N(p)=\frac{1}{\sqrt{r^{2}+ r^{2}(f'(r))^{2}}} (-f'(r) x,-f'(r) y,r ),
\end{align}
and  euclidean principal curvatures and euclidean mean curvature
\begin{align}\label{principal-curvatures}
k_{1}(p,E)=-\frac{f''(r)}{(1+ (f'(r))^{2})^{\frac{3}{2}}}, \qquad k_{2}(p,E)=-\frac{f'(r)}{r(1+ (f'(r))^{2})^{\frac{1}{2}}}, \quad H(p,E)=k_{1}(p)+ k_{2}(p).
\end{align}
\end{remark}

\subsection{Horizontal mean convex  sets} 
We introduce the definition of horizontal mean convex set and horizontal strictly mean convex set.
\begin{definition}\label{meanconvex} We say that a $C^2$ set $B$  is horizontal mean convex with constant $\delta\geq 0$,  if   $H_{\mathbb{H}}(p,B)\geq\delta$ for every $p\in \partial B\setminus (\partial B)_0$, where $(\partial B)_0$ denotes the set of singular points of the boundary. 

If $\delta>0$, then $B$ is called horizontal  strictly mean convex with constant $\delta$.  

\end{definition} 
\begin{remark}\upshape \label{remconvex}  
If  $E$ is  an axially symmetric set (as described in Remark \ref{remarkaxi}) and it is (euclidean) convex, then  it is horizontal mean convex. Indeed,   since $f''\leq 0$ and $f'\leq 0$,  we deduce using \eqref{radialmc} that $H_{\mathbb{H}}(p, E )\geq 0$. 
\end{remark} 
In the following we collect  some examples of horizontal mean convex and strictly mean convex sets. 

\begin{remark}\upshape[Hyperplanes]\label{planecurvature}
Recalling the formulas for the horizontal normal provided in Remark~\ref{planenormal}
 we apply  \eqref{divnormal2} and \eqref{horizontalmc} away from $\Sigma_0$ to show that every plane has zero  horizontal mean curvature.
Every horizontal plane $\{p\ |\ z=k\}$ has zero horizontal mean curvature since 
\[H_{\mathbb{H}}(p,  \{ z\leq k\} )=\partial_x\left(\frac{-y}{\sqrt{x^2+y^2}}\right)+\partial_y\left(\frac{x}{\sqrt{x^2+y^2}}\right)=0. \] The same argument gives that the planes  $\{p\ |\ ax+by+z=k\}$ has zero mean curvature since
\[ \partial_x\left(\frac{a-y/2}{\sqrt{(a-y/2)^2+(b+x/2)^2}}\right)+\partial_y\left(\frac{b+x/2}{\sqrt{(a-y/2)^2+(b+x/2)^2}}\right)=0.\]
The vertical planes $\{p\ |\ ax+by=k\}$ have constant horizontal normal, so they have trivially zero  horizontal mean curvature.
 
 Therefore all  hyperplanes are horizontal mean convex sets. 

\end{remark}

\begin{remark}\upshape[Balls]\label{spherecurvature}
We compute  the horizontal mean curvature of some different  balls, based on the previous formula \eqref{radialmc}, see \cite[Section 4.3]{CapognaDanielli}. 

\begin{enumerate}\item
Let us consider the boundary of euclidean ball $  B_E(0,R)=\{(x,y,z)\, | \,  x^2+y^2+z^2< R^2\}$. 
Then 
\[H_{\mathbb{H}}(p, B_E(0,R)) = \frac{2(4+R^2)}{\sqrt{x^2+y^2}(4+z^2)^{3/2}}.
\] Note that the horizontal curvature explodes at the singular points $(0,0,\pm R)$. In particular a euclidean ball is a horizontal strictly mean convex with constant $\frac{2 }{R\sqrt{ 4+R^2}}$. 
\item The Kor\'anyi ball  $ B_K(0,R)=\{(x,y,z) \, | \,  (x^2+y^2)^2+16 z^2<R^4\}$ has curvature at $p\in \partial B_K(0,R)\setminus \{(0,0, \pm R^2/4)$ given by 
\begin{equation}\label{curvaturekoriani}H_{\mathbb{H}}(p, B_K(0,R)) =\frac{3\sqrt{x^2+y^2}}{R^2}.
\end{equation} 
The Kor\'anyi ball is  convex and moreover it is horizontal mean convex (but not strictly mean convex) since $0\leq H_{\mathbb{H}}(p, B_K(0,R))\leq \frac{3}{R}$ for all $p\in\partial B_K(0,R)\setminus (\partial B_K(0,R))_0$.

Note that the euclidean principal curvatures (see \eqref{principal-curvatures}) are given by\begin{align*}
k_{1}(p, B_K(0,R))= \frac{4r^2 (3R^4-r^4)}{(4(R^{4}-r^{4})+r^{6})^{\frac{3}{2}}}, \qquad  k_{2}(p, B_K(0,R))= \frac{r^2}{(4(R^{4}-r^{4})+r^{6})^{\frac{1}{2}}}
\end{align*} where $r=\sqrt{x^2+y^2}$. 
The  euclidean principal curvatures vanish for $r \to 0$.

\item The boundary of the Carnot--Carathéodory ball  centered at $0$, that is $\partial B_{CC}(0,R)$, can be parametrized  recalling that  a length minimizing curve between $0$ and $(x,y,z)$   is the lift of a circular arc joining the origin in $\R^2$ with $(x,y)$, whose convex hull has area $z$ (see Remark~\ref{geodesic}).
The parametrization   is given by
\[\sqrt{x^2+y^2}=\frac{2}{c}\sin(cR/2) \qquad z= \frac{cR-\sin cR}{2c^2}.\]
 The points with maximal and minimal height (that is the points such that $\frac{dz}{dc}(c)=0$), correspond to $cR=\pm \pi$ and have coordinates $(x,y,z)$ with $\sqrt{x^2+y^2}=\frac{2R}{\pi}<R$, and $z=\pm \frac{R^2}{2\pi}$. 
The points on the plane $\{z=0\}$ corresponds to $c=0$, and satisfy $\sqrt{x^2+y^2}=R$, whereas the points on the $z$-axis  correspond to $c=\pm (2\pi)/R$ and have coordinate $(0,0, \pm R^2/(4\pi))$.  It is possible to compute explicitly 
\begin{equation}\label{curvaturecc}H_{\mathbb{H}}(p, B_{CC}(0,R))=\frac12 \frac{c}{2\sin(cR/2)}\frac{\sin cR-cR\cos cR}{\sin (cR/2)-cR\cos (cR/2)/2}. 
\end{equation}
Note that if  $p=(x,y,0)\in \partial B_{CC}(0,R)$, then $c=0$ and  $H_{\mathbb{H}}(p, B_{CC}(0,R))=\frac{4}{R}$. At points of maximal  and minimal height, $cR=\pm \pi$ and $H_{\mathbb{H}}(p, B_{CC}(0,R))=\frac{\pi^2}{4R}$.
On the other hand, if $p=(x,y,z)\in \partial B_{CC}(0,R)$ satisfies  $\sqrt{x^2+y^2}\to 0$, then  $c\to\left( \frac{2\pi}{R}\right)^{-}$ and we obtain  $H_{\mathbb{H}}(p, B_{CC}(0,R))\to -\infty$. 

Thus, the Carnot--Carathéodory ball is  not horizontal mean convex. 
\end{enumerate}
\end{remark}
\subsubsection{Pansu ball}  \label{Pansuspheredef}
The Pansu sphere, as introduced in Definition \ref{defpansu}   is the surface obtained by rotating around the $z$-axis  the geodesics connecting $0$ with $(0,0,\pi R^2)$. In order to center the sphere in the origin, the set is vertically translated and dilated.
So, we will consider the Pansu sphere of parameter $R$ as the surface obtained by rotating around the $z$-axis  the geodesics connecting $(0,0, - \pi R^2/8)$ with $(0,0, \pi R^2/8)$, which coincides with the singular points of the sphere. 
By recalling the description of geodesic in Remark \ref{geodesic}, it is possible to show (see e.g. \cite{r}) that the Pansu ball is defined as 
\begin{equation}B_P^R=\{(x,y,z)\ | \  |z |<f_R(\sqrt{x^2+y^2}) \}\end{equation} where 
\[f_R(r)= \frac14 R^2\arccos \frac{r}{R}+\frac14 r\sqrt{R^2-r^2}= \frac14 R^2\left(\arccos \frac{r}{R}+  \frac{r}{R} \sqrt{1-\frac{r^2}{R^2}}\right).\]

The horizontal mean curvature of  $\partial B_P$ is constant: 
\begin{equation}\label{curvaturepansui}H_{\mathbb{H}}(p,  B_P^R)= \frac{2}{R}.
\end{equation} Note that $\delta_{1/R}B_P^R= B_P^1$ or equivalently $B_P^R=\delta_{R} B_P^1$ and the boundary of the Pansu ball is $C^{2}$ but not $C^{3}$ (\cite{CapognaDanielli}).
The Pansu ball is a convex set which is horizontal   strictly mean convex.

Observe that for the Pansu ball the euclidean principal curvatures can be computed using \eqref{principal-curvatures} obtaining  \begin{align*}
k_{1}(p,B^{R}_{P})= \frac{4r (2R^{2}-r^{2})}{(4(R^{2}-r^{2})+r^{4})^{\frac{3}{2}}}, \qquad 
k_{2}(p,B^{R}_{P})=\frac{r }{ {(4(R^{2}-r^{2})+r^{4})^{\frac{1}{2}}}}
\end{align*}  where $r=\sqrt{x^2+y^2}$.
In particular  the euclidean  principal curvatures vanish for $r \to 0$.

In \cite[Theorem 6.1]{rr} an Alexandrov theorem in $\mathbb{H}$ is proved about the Pansu spheres:
they   are the only $C^2$ immersed compact connected surfaces  with positive  constant horizontal mean curvature and the  regular part of every Pansu sphere  is foliated by horizontal geodesics of curvature equal to the constant curvature of the sphere.

\begin{remark}\label{remarkbb} \upshape Observe that by the Ball-box theorem there exists $\lambda_1<1$, $\lambda_2\geq 1$ such that $B_P^{\lambda_1R}\subseteq B_{CC}(0, R)\subseteq B_P^{\lambda_2R}$.
\end{remark}

 \subsection{Horizontal perimeter and functions of bounded variations}

 It is possible to adapt the definition of  functions of bounded variation  to the case of $\mathbb{H}$, see \cite{serra}, \cite[Chapter 5]{CapognaDanielli}. 
 We define the notion of horizontal variation for a measurable function $f:\mathbb{H}\to \R$:
 \begin{eqnarray*}
\Var_{\mathbb{H}}(f;\Omega)&:=
& 
\sup\left\{ \int_{\Omega} f \Div_{\mathbb H}\phi dp\ | \  \phi \in C^{1}_c(\Omega, \R^2), \sqrt{(\phi_1)^2+(\phi_2)^2}\leq 1\right\}\\ 
&=& 
\sup\left\{ \int_{\Omega} f \Div_g\phi_{\mathbb H} dp\ | \  \phi \in C^{1}_c(\Omega, \R^3), |\phi_\mathbb{H}|_g\leq 1\right\},
\end{eqnarray*}
where $\Omega\subseteq \mathbb{H}$ denotes an open set.
Recall that for $\phi \in C^{1}(\Omega, \R^2)$, it holds  $\Div_{\mathbb H}\phi =X(\phi_1)+Y(\phi_2)$, and for $\phi \in C^{1}(\Omega, \R^3)$ it holds $ \Div_g\phi_{\mathbb H} =\Div\phi_{\mathbb H} $ (cp. with \eqref{div-equi}).

\begin{definition}
Let $\Omega\subseteq \mathbb{H}$ be an open set and $f: \Omega \subset \mathbb{H} \rightarrow \mathbb{R}$. Then $f \in BV_{\mathbb{H}}(\Omega)$   if $f \in L^{1}(\Omega)$ and if
 $\Var_{\mathbb{H}}(f;\Omega)<+\infty$.

Let  $\mathcal{M}$  be  the class of all measurable sets in $ \mathbb{H}$.  A set $E \in\mathcal{M}$  has  finite $\mathbb{H}$ perimeter in $\Omega$ if $\Per(E, \Omega):=\Var_{\mathbb{H}}(\chi_E;\Omega)<+\infty$, where $\chi_E$ is the characteristic function. More explicitly:
\begin{align}\label{perimeter}
P_{\mathbb{H}} (E;\Omega) =  \sup_{\phi \in C^{1}_c(\Omega,\R^2), |\phi|\leq 1} \int_{E}  \Div_{\mathbb{H}}\phi\; dp.\end{align}
If $\Omega=\mathbb{H}$ we denote $P_{\mathbb{H}} (E;\mathbb H)=P_{\mathbb{H}} (E)$. 
Moreover $E$ is a $\mathbb{H}$-Caccioppoli set in $\Omega$ if  $\chi_{E} \in BV(\mathbb{H})(\Omega)$. We will denote 
\begin{align}\label{mengeP} \mathcal{P}:=  \{ F \in \mathcal{M} \, : \, \chi_F\in BV(\mathbb{H}) \}.
\end{align}
\end{definition}

The perimeter functional $P_{\mathbb{H}}$ is lower semicontinuous with respect to the $L^1_{loc}$ convergence (\cite[Proposition 5.6]{CapognaDanielli}). 

\begin{remark} If $F$ is a set of finite perimeter then $\mathcal{H}_{CC}^{3} (\partial^{e}F \setminus F(1/2))=0$, where  $\mathcal{H}^{s}_{CC}$ denotes the $s$-dimensional Hausdorff measure obtained from $d_{CC}$ (see \cite[Proposition 1.23]{lr}).
\\ If $E$, $E'$ are  measurable sets, then   
\begin{align} \label{eq30} P_{\mathbb{H}} (E, \Omega) &=P_{\mathbb{H}} (\mathbb{H}\setminus E, \Omega),\\ \label{eq31}
P_{\mathbb{H}} (E, \Omega) &= P_{\mathbb{H}} (E', \Omega) \qquad \text{ whenever } |(E \triangle E') \cap \Omega| =0
\end{align} 
where $E \triangle E'= (E \setminus E') \cup (E' \setminus E)$ is the  symmetric difference between $E$ and $E'$.

The sets $E$ and $E'$ are said to be equivalent if  $|E \triangle E'|=0$.
\end{remark}

For regular sets, the horizontal perimeter has an   explicit integral representation in terms of the underlying Riemannian  geometry.  
We recall the divergence theorem for the Riemannian space $(\mathbb H, g)$. Let $M\subseteq \mathbb{H}$ a $C^1$ submanifold with boundary, and $U$ a $C^1$ vector field then
\begin{equation}\label{divergencetheorem}\int_M \Div_g U dp= \int_{\partial M} U\cdot_g N_g d\sigma, \end{equation}
where $N_g$ is the Riemannian outer normal at $\partial M$ and $d\sigma$ is the Riemannian area element on $\partial M$.

\begin{proposition}[Area of a surface] 
If $E$ is set of class $C^1$, then
\begin{equation}\label{area}
\Per(E;\Omega) =\sup_{\phi \in C^{1}_c(\Omega,\R^2), |\phi|\leq 1} \int_{\partial E}  \phi\cdot_g N_g   d\sigma=  \int_{\partial E\cap \Omega} \nu_{\mathbb H}\cdot_g N_g  d\sigma= \int_{\partial E\cap \Omega} |N_{\mathbb H}|_g  d\sigma 
\end{equation} 
where $N_{\mathbb{H}}$ is defined in \eqref{nh} and $d\sigma$ is the   Riemannian area element on $\partial E$. 

In terms of Euclidean unit normal $N=\frac{\nabla f}{|\nabla f|}$ and standard 2-dimensional Hausdorff measure $\mathcal{H}^{2}$ we have
\begin{align}\label{area2}
\Per(E;\Omega)
=\int_{\partial E\cap \Omega}  \nu_{\mathbb{H}} \cdot N(p)  d\mathcal{H}^{2}(p) =
\int_{\partial E\cap \Omega} \frac{\sqrt{ (Xf)^{2}+ (Yf)^{2}}}{|\nabla f|}   d\mathcal{H}^{2} .
\end{align}

\end{proposition} 
For the proof we refer to \cite[Proposition5.7]{CapognaDanielli} or \cite[Proposition 3.1]{Monti}.

\begin{remark}[Invariance properties of the perimeter] \label{inv} \upshape The  horizontal perimeter by its very definition  is invariant with respect to any isometry of $\mathbb{H}$ which  leaves invariant the horizontal distribution. So  in particular it is invariant with respect to left translations that is \[P_{\mathbb H} (E;\Omega)=P_{\mathbb H} (L_{p}(E); L_{p} (\Omega))\]  and rotations around the $z$-axis . 
Indeed since $X, Y$ are left invariant, then $(\Div_{\mathbb{H}} \phi)\circ L_p= (\Div_{\mathbb{H}} \phi\circ L_p)$. 
Moreover,  let  $\theta \in (0, 2\pi)$ and $R_\theta$ the associated rotation around the $z$-axis   of angle $\theta$, 
then \[ P_{\mathbb{H}}(E)= P_{\mathbb{H}}(R_{\theta}E). \]
\end{remark}

\begin{remark}\upshape  \label{dilatations} The horizontal perimeter  is    homogeneous of degree $3$ with respect to dilations of sets, see \cite[Lemma 5.5]{CapognaDanielli} that is for all   $\lambda>0$, 
\[  P_{\mathbb H} (\delta_\lambda E;\delta_\lambda\Omega)=\lambda^3 P_{\mathbb H} ( E;  \Omega).\] \end{remark} 
The horizontal perimeter is also subadditive, see \cite{CapognaDanielli}. The proof of this fact follows the same arguments as for the standard perimeter. 
 \begin{lemma} For $E,F \in \mathcal{P}$ we have 
\begin{align}\label{propPer}
P_{\mathbb{H}} (E\cap F) + P_{\mathbb{H}} (E\cup F) \leq P_{\mathbb{H}} (E) +
P_{\mathbb{H}} (F).
\end{align}
\end{lemma}

Also in the Heisenberg group we may formulate the classical isoperimetric problem:
\begin{equation}\label{iso}\min_{E\in  \mathcal{M}, |E|=1}P_{\mathbb{H}}(E).
\end{equation} 
Existence of isoperimetric sets is proved in \cite{lr} (see also \cite{flm}, \cite{Monti}).   

\begin{theorem}[Isoperimetric inequality]
If $F$ is any set with finite $\mathbb{H}$-perimeter  then  either $F$ or $\mathbb{H}\setminus F$ is a $\mathbb{H}$-Caccioppoli set  and there exists a universal constant $C>0$ such that \begin{equation}\label{isoineq}
 P_{\mathbb{H}}(F)
\geq C  \min ( |F|, |\mathbb{H}\setminus F|)^{3/4}.\end{equation}
\end{theorem}
\begin{proof} See \cite[Theorem~2.2]{lr}.
\end{proof}
\begin{remark}\label{pansuconstant} \upshape  The Pansu sphere $B_P^R$ as introduced in Definition \ref{defpansu} is conjectured  to be the solution to the  isoperimetric problem \cite{p}. 
We have that 
\[|B_P^R|= \frac{3}{16}\pi^2 R^4\qquad \Per(B_P^R)=\frac12 \pi^2 R^3\] so the isoperimetric constant is conjectured to be
\[C:=\frac{4\sqrt{\pi}}{3^{3/4}}.\]
The conjecture  was proved in the class of axially symmetric sets in \cite{m}, in the convex setting in \cite{mr},  in the smooth setting \cite{rr}, and in the axially symmetric cylinders \cite{r}.  
\end{remark} 
\begin{theorem}[CC-Balls are domains of isoperimetry]\label{teoBallIso}
There exists a constant $C>0$ such that for any measurable  set $E$  and any ball $B_{CC}(p, r)$ (as in \eqref{ccball}) such that
$\emptyset \neq E\cap B_{CC}(p, r)$ there holds
$$ P_{\mathbb{H}}(E\cap B_{CC}(p, r),B_{CC}(p, r))\geq C \min \{ |E\cap  B_{CC}(p, r)|, |B_{CC}(p, r)\setminus E| \}^{3/4}. $$
\end{theorem}
\begin{proof}
See \cite[Theorem 2.15 and Definition 2.13]{lr}.
\end{proof}

Finally let us recall an  important compactness theorem.
 \begin{theorem}[Compactness in $BV_{\mathbb{H}}$]\label{CompactnessStrong}
 Let $(E_{i})$ be a sequence of measurable set with $\sup_{i}\Per(E_{i}) < \infty.$
 Then, up to passing to a subsequence,  $(E_{i})$ converges in $L^{1}_{loc}(\mathbb{H})$ to a measurable set $E_{\infty}$ with finite $\mathbb{H}$-perimeter.
 
 Moreover, if there exists a bounded open set  $\Omega$  such that $E_i\subseteq \Omega$ for all $i$, then $E_i$ converges up to subsequences in $L^{1}(\mathbb{H})$ to  $ E_\infty\subseteq  \overline{\Omega}$  with finite $\mathbb{H}$-perimeter.
 \end{theorem}
 \begin{proof} See \cite[Theorem~8.4]{CapognaDanielli}. \end{proof}

\subsection{Outward minimizing sets}
In this section we introduce the notion of outward minimizing sets and we provide some relevant examples.  Upon recalling \eqref{mengeP} we give the following
\begin{definition}\label{def1}
Let $\Omega \subset \mathbb{H}^{1}$ be an open set. We say that a set $E \in \mathcal{P}$, $E \subset \Omega$, is \emph{outward minimizing} in $\Omega$ if
\begin{align}
P_{\mathbb{H}} (E) \leq P_{\mathbb{H}} (F) \qquad \text{ for all  $F \in \mathcal{P}$ with } E \subset F \subset \Omega.
\end{align}
\end{definition}
An equivalent characterization of the outward minimality is the following:  
\begin{lemma}\label{lemma1.1}
The set $E \in \mathcal{P}$ is outward minimizing in $\Omega$ if and only if
\begin{equation}\label{char}
P_{\mathbb{H}} (E\cap G) \leq P_{\mathbb{H}} (G) \qquad  \text{ for all } G \subset \Omega, G \in \mathcal{P}.
\end{equation} 
Moreover if $E_n$ is a sequence of  outward minimizing sets in $\Omega$, and $E_n\to E$ in $L^1(\Omega)$, then $E$ is outward minimizing in $\Omega$. 
\end{lemma}
\begin{proof}
One can argue exactly as in  \cite[Lemma~2.5]{DePhilLaux} using \eqref{propPer}. For sake of completeness we report it here. 
If condition \eqref{char} holds, then we apply it to $G\supseteq E$ and get immediately $P_{\mathbb{H}} (E) \leq P_{\mathbb{H}} (G)$. On the other hand, assume $E$ outward minimizing. 
Since $E\cup G\supseteq E$ by outward minimality $P_{\mathbb{H}} (E) \leq P_{\mathbb{H}} (E\cup G)$.
By \eqref{propPer}, $P_{\mathbb{H}} (E\cap G) +P_{\mathbb{H}} (E\cup G) \leq P_{\mathbb{H}} (G)+ P_{\mathbb{H}} (E)$.  So, putting together the two inequalities we get 
 $ P_{\mathbb{H}} (E\cap G) +P_{\mathbb{H}} (E ) \leq  P_{\mathbb{H}} (E\cap G) +P_{\mathbb{H}} (E\cup G) \leq P_{\mathbb{H}} (G)+ P_{\mathbb{H}} (E)$, and so 
\eqref{char}. 
 
The stability under $L^1$ convergence is a consequence of the characterization of outward minimality given in \eqref{char} and of the lower semicontinuity of $\Per$ with respect to $L^1$ convergence, see \cite[Proposition 2.7]{cn}.
\end{proof}

We recall that outward minimizing sets satisfy an upper density estimate.
 
\begin{lemma}[Upper density estimate]\label{upperdensityest}
 Let $E \in \mathcal{P}$ be an outward minimizing set in $\mathbb{H}$.
 There exist  a constant $ \theta >0$  such that for all $p \in \mathbb{H}$, $r>0$ there holds 
  $$  P_{\mathbb{H}}(E, B_{CC}(p,r)) \leq  \theta r^{3} $$
where  $B_{CC}(p,r)$  is defined in \eqref{ccball}. 
 \end{lemma}
 \begin{proof}The proof  follows exactly as in \cite[Lemma~4.1]{lr} by substituting the property that $E$ is an isoperimetric set  with the property that $P_{\mathbb{H}}(E) \leq P_{\mathbb{H}}(E \cup B_{CC}(p,r))$ for all $r>0$.
 \end{proof}

We provide next some relevant  examples of outward minimizing sets in the Heisenberg group, using a calibration argument (see \cite{r}, and \cite{rr}).

First of all we recall the definition of area-minimizing surface.
\begin{definition} A surface $\Sigma$ is area-minimizing if any region $M\subseteq\Sigma$ satisfies $P_{\mathbb{H}} (M)\leq P_{\mathbb{H}} (M')$ for  any other $C^1$ compact surface $M'$ with $\partial M=\partial M'$.\end{definition} 
In \cite[Theorem~4.17]{rr} it is proved that if an immersed surface is area-minimizing  then its horizontal 
mean curvature is zero (or constant under a volume constraint) and the characteristic curves (with are the integral curves of $J(\nu_{\mathbb{H}})$) meet orthogonally the singular curves (if there are any). 
Moreover  it is proved that every non vertical plane in $\mathbb{H}$ is globally area minimizing, whereas in \cite{dgar}  vertical Euclidean
planes are characterized as the only $C^2$ area-minimizing entire graphs in $\mathbb{H}$  with empty singular set.     
In particular, by arguing first for smooth sets and then proceeding by regularization, this implies  that for every set $E\in \mathcal{P}$ with finite $\mathbb{H}$ perimeter and for every half-space $\Pi$  it holds
\[\Per(E\cap \Pi)\leq \Per(E).\] 
%


\subsubsection{A family of  outward minimizing sets} 


We introduce a class of smooth sets, and we show, by using a calibration type argument, that they are outward minimizers. 
 \begin{definition}\label{definitionb} We denote by $\mathcal{B}$ the class of  open bounded sets $B\subseteq \mathbb H$ of class $C^2$  which are  axially symmetric (with respect to the z-axis),   strictly starshaped w.r.t. $0$ according to Definition \ref{strictlystarshaped}, horizontal   mean convex according to Definition \ref{meanconvex},  and moreover satisfy the following   regularity condition:
 there exists $H_0=H_{0}(B)>0$ such that \[H_{\mathbb H}(p, B)\leq H_0\qquad \forall\ p\in \partial B\setminus (\partial B)_0.\]  
 \end{definition} 
The most important example  in the class $\mathcal{B}$ is the Pansu ball of parameter $R>0$, centered at the origin, described in Section \ref{Pansuspheredef}.   Moreover by  Remark~\ref{spherecurvature},  
K\`oranyi balls with radius $R$ are in the class $\mathcal{B}$.  
On the  other hand the Euclidean ball of radius $R$ is strictly starshaped with respect to $0$ and horizontal strictly mean convex with constant $\delta=\frac{2}{R\sqrt{4+R^2}}$, but does not satisfy the upper bound on the curvature.

\begin{proposition}\label{convexset} Let $B\in\mathcal{B}$,  horizontal mean convex  with constant $\delta\geq 0$.  Let   $B^\lambda:= \delta_\lambda B$, for $\lambda>0$. 
Then  $B$ is   an  outward minimizing set in $\mathbb{H}$, that is for any $F\in \mathcal{P}$ such that $B\subseteq F$ it holds 
\begin{equation}\label{outwardB} \Per(F)\geq \Per(B)+\frac{\delta}{\lambda_F}|F\setminus B|\end{equation}  where      $\lambda_F:=  \inf\{\lambda \,|\, F\subseteq B^\lambda\}\geq 1$, with $\lambda_F=+\infty$ if $F$ is unbounded.  Moreover 
    for general $F\subseteq \mathbb{H}$ \begin{equation}\label{outwardB2} \Per(F)\geq \Per(F\cap B)+\frac{\delta}{ \max\{\lambda_F,1\}}|F\setminus B|\end{equation} 
   \end{proposition}
 \begin{remark}\upshape  Proposition \ref{convexset} applies to the Pansu ball $B_P^R$ with $\delta=\frac{2}{R}$ and to the 
K\`oranyi ball with radius $R$ with $\delta=0$.  

\end{remark}

\begin{proof}    
First of all we observe that if  $B$ is $C^2$ and axially symmetric, then the singular set $(\partial B)_0$ is contained in the $z$-axis as proved in \cite[Section 4]{rr1}. 
We divide the proof of the first statement in several steps. 

\noindent {\bf Step 1: definition of the homogeneous   gauge.}
Since $B$ is  strictly starshaped w.r.t. $0$ as in Definition \ref{strictlystarshaped} then the family of  dilations $\delta_\lambda B:=B^\lambda$ for $\lambda\geq 0$ produce a foliation of $\mathbb H$ as proved in Lemma \ref{lemmastellati}. 
We may associate to $B$  an  homogeneous  gauge: 
for $p\ne0$ let us define 
\begin{equation}\label{gaugeB}\rho_B(p) \text{ be the unique positive number such that }
p\in\partial B^{\rho_B(p)},
\end{equation} 
and set $\rho_B(0)=0$. Thus $B^\lambda =\{\rho_B<\lambda\}$ and $\rho_B(\delta_\lambda p)=\lambda\rho_B(p)$. The function $\rho_B$ is a continuous strictly positive function, homogeneous of degree $1$ with respect to dilations. 
 By compactness of the unit euclidean sphere  we have that the map $p\mapsto \frac{1}{\rho_B(p)}$ is locally integrable at $0$, by passing to spherical coordinates and using the fact that the 
homogeneous dimension in $\mathbb H$ is  $4$, so greater than $1$.  

\noindent {\bf Step 2: extension of the outer horizontal unit normal.}
Away from the $z$-axis define the horizontal vector field
\begin{equation}\label{horizontalnormalestesa}
U_B(p)=\nu_{\mathbb H}^{B^{\rho_B(p)}}(p)
\end{equation}  where $\nu_{\mathbb H}^{B^{\rho_B(p)}}(q)$ is the horizontal exterior  normal to $ B^{\rho_B(p)}$ at $q\in \partial  B^{\rho_B(p)}$. 
 An explicit expression for $U_{B}$ can be computed using \eqref{radialnormal}. In particular this vector field is $C^{1}$ away from the $z$-axis.
Therefore, away from the $z$-axis we have, recalling  \eqref{horizontalmc} and Remark~\ref{remarkhomoH},  \begin{equation}\label{divU}
|U_B|_g=1,\qquad \frac{\delta}{\rho_B(p)}\leq  \operatorname{div}_{\mathbb H}U_B(p)=H_{\mathbb H} (p, B^{\rho_B(p)})\leq \frac{H_0}{\rho_B(p)}.
\end{equation}
We observe  that $p\mapsto \operatorname{div}_{\mathbb H}U(p) $ is locally integrable at $0$.  

\noindent {\bf Step 3: calibration argument:
for every $G\in\mathcal P$  there holds 
\begin{equation}\label{pansu-calibration-ineqB}
P_{\mathbb H}(G)\ge \int_G \operatorname{div}_{\mathbb H}U_B(p)dp\geq  \int_G\frac{\delta}{\rho_B(p)}\,dp.
\end{equation}
Moreover equality holds when $G=B^\lambda$, $\lambda>0$:
\begin{equation}\label{pansu-calibration-eqB}
P_{\mathbb H}(B^\lambda)=\int_{B^\lambda} \operatorname{div}_{\mathbb H}U_B(p)dp.
\end{equation}}

The calibration argument is based  on \cite[Theorem 3.1]{r}. We consider two family of  smooth functions  $\phi_\epsilon(p), \psi_\epsilon(p):\mathbb{H}\to[0,1]$,  defined as follows. 

Let us fix a smooth function $\phi_\eps:\R^2\to [0,1]$ such that 
\begin{eqnarray*} 
\phi_\epsilon(x,y)\equiv 0 & \text{ for  }\sqrt{x^2+y^2}\leq \epsilon^2 \text{ and } \sqrt{x^2+y^2}\geq \frac{1}{\sqrt{\epsilon}}+1, \\
\phi_\epsilon(x,y)\equiv 1  & \text{ for  } \epsilon \leq  \sqrt{x^2+y^2}\leq \frac{1}{\sqrt{\epsilon}},\\
 |\nabla \phi_\epsilon |\leq 2/\epsilon &  \text{ for  }\epsilon^2\leq  \sqrt{x^2+y^2}\leq \epsilon,\\
  |\nabla \phi_\epsilon |\leq 2 &  \text{ elsewhere.   }\end{eqnarray*} 
It is straightforward to observe that we may extend $\phi_\eps$ to be a smooth function on $\mathbb{H}$, by setting $(x,y,z)\in \mathbb{H}\mapsto \phi_\eps(x,y)$. Moreover $\nabla\phi_\eps=\nabla_{\mathbb H}\phi_\eps$ and $|\nabla_{\mathbb{H}}\phi_\eps|=|\nabla \phi_\eps|_g$.

Now we consider a smooth function $\psi_\eps:\R\to[0,1]$  such that 
\[  \psi_\epsilon(z)\equiv 1, \text{ for  } |z|\leq\frac{1}{\sqrt{\epsilon}}  \quad  \psi_\epsilon(z)\equiv 0, \text{ for  } |z|\geq\frac{1}{\sqrt{\epsilon}}+1   \] with   $| \psi_\epsilon'|\leq 2$.  Also in this case we may extend  $\psi_\eps$ to be a smooth function on $\mathbb{H}$, by setting $(x,y,z)\in \mathbb{H}\mapsto \psi_\eps(z)$. Moreover $|\nabla_{\mathbb H}\psi_\eps|_g= |\psi_\eps'| \sqrt{\frac{x^2+y^2}{4}}$.

Therefore, the vector field $U_B\phi_\epsilon\psi_\epsilon$ is still horizontal, and it is $C^1$,  compactly supported,  identically $0$ on the singular set $\R(0,0,1)$.  Therefore by definition of perimeter 
\eqref{perimeter} we obtain that for every  $G\in\mathcal P$: 
\[\Per (G)\geq \int_G \Div_{\mathbb H}(U_B\phi_\epsilon\psi_\epsilon)(p) dp.\] 

At any point $p=(x,y,z)$ with $z\neq 0$   we have 
\begin{eqnarray} \label{divUeps}   \Div_{\mathbb H}(U_B\phi_\epsilon\psi_\epsilon)(p) &\\ =&  \phi_\eps \psi_\eps \Div_\mathbb{H}U_B+ \psi_\eps [(U_B)_1\partial_x\phi_\eps +(U_B)_2\partial_y\phi_\eps]+\phi_\eps\left[-\frac{y}{2}(U_B)_1 +\frac{x}{2} (U_B)_2 \right]\psi_\eps'(z).\notag    \end{eqnarray}

 First of all we assume  $G\subseteq \mathcal{C}_s=\{(x,y,z) \,|\, x^2+y^2\leq s^2, |z|\leq s\}$ for some $s>0$. 
We claim that 
\begin{equation}\label{claim1}  \lim_{\epsilon\to 0}\int_{G}  \Div_{\mathbb H} (U_B \phi_\epsilon\psi_\epsilon ) dp=\lim_{\epsilon\to 0}\int_{G} \phi_\epsilon\psi_\epsilon \Div_{\mathbb{H}} U_B dp.\end{equation} 
Note that if the claim holds true, we infer \eqref{pansu-calibration-ineqB}. Indeed $\Per(G)\geq  \lim_{\epsilon\to 0}\int_{G}  \Div_{\mathbb H} (U_B \phi_\epsilon\psi_\epsilon ) dp$ and moreover
$\lim_{\epsilon\to 0}\int_{G} \phi_\epsilon\psi_\epsilon \Div_{\mathbb{H}} U_B dp= \int_{G}   \Div_{\mathbb{H}} U_B dp$
by the dominated convergence theorem since $\chi_G  \Div_{\mathbb{H}} U_B\in L^{1}$.

To prove the claim \eqref{claim1} we recall the formula \eqref{divUeps}. 
First of all we observe that
\begin{equation}\label{uno}\lim_{\eps\to 0} \int_G \phi_\eps\left[-\frac{y}{2}(U_B)_1 +\frac{x}{2} (U_B)_2 \right]\psi_\eps'(z) dp=0 \end{equation}  due to the dominated convergence theorem, since $\chi_G\in L^1$,  
$\phi_\eps\left[-\frac{y}{2}(U_B)_1 +\frac{x}{2} (U_B)_2 \right]\psi_\eps'(t)\chi_G$ is uniformly bounded since $G\subseteq \mathcal{C}_s$   and is converging to $0$. 
Next, we have  that  \begin{equation}\label{due} \left|\int_{G} [(U_B)_1\partial_x\phi_\eps+(U_B)_2\partial_y\phi_\eps]\psi_\epsilon  dp \right| \leq \int_{G}  \left|\nabla(\phi_\epsilon)\right| \psi_\epsilon dp   \leq \int_{G} \chi_{A_\eps} \frac{2}{\epsilon}\psi_\epsilon dp  +  \int_{G} \chi_{B_\eps} 2\psi_\epsilon
  dp\end{equation}  where $A_\eps=\{p \ | \epsilon^2\leq \sqrt{x^2+y^2} \leq \epsilon \}$ and $B_\eps=\{p \ | \frac{1}{\sqrt{\epsilon}}\leq \sqrt{x^2+y^2} \leq \frac{1}{\sqrt{\epsilon}}+1\}$.   
  Since $2\chi_{B_\eps}\psi_\eps\to 0$ and $\chi_G\in L^1$, we infer by dominated convergence that 
\begin{equation}\label{tre} \lim_{\eps\to 0}  \int_{G} \chi_{B_\eps} 2\psi_\epsilon
  dp=0.\end{equation}  On the other hand passing to cylindrical coordinates 
\begin{equation}\label{quattro} \lim_{\eps\to 0}  \int_{G} \chi_{A_\eps} \frac{2}{\epsilon}\psi_\epsilon dp\leq \int_{\eps^2}^{\eps} 2\pi  r \int_{-\frac{1}{\sqrt{\eps}}-1}^{\frac{1}{\sqrt{\eps}}+1} \frac{2}{\eps} dzdr = \lim_{\eps\to 0} \int_{\eps^2}^{\eps} r \frac{8\pi}{\eps }\left(\frac{1}{\sqrt{\eps}}+1\right)dr=0. \end{equation}
 
Therefore, by \eqref{divUeps}, \eqref{uno},\eqref{due}, \eqref{tre},  \eqref{quattro} the claim  \eqref{claim1} and thus \eqref{pansu-calibration-ineqB} is proved for $G$ bounded.  
If $G$ is not bounded, observe that by coarea formula  \cite[Theorem 6.11]{CapognaDanielli}  there exists $s_k\to +\infty$ such that $P_{\mathbb H}(G_k)\to P_{\mathbb H}(G)$, where $G_k=G\cap B_{CC}(0,s_k)$ (an outline of this step is given for instance in \cite[Remark 3.12]{Maggi}). 
Since $\lim_{s_k} \int_{G_k}  \operatorname{div}_{\mathbb H}U_B(p)dp \geq \int_G \operatorname{div}_{\mathbb H}U_B(p)dp $ by Fatou's lemma,   and 
 $P_{\mathbb H}(G_k)\to P_{\mathbb H}(G)$, we conclude  \eqref{pansu-calibration-ineqB} also for $G$. 

On the other hand,   equality \eqref{pansu-calibration-eqB}  is obtained by the Gauss Green formula  for $C^2$ sets \eqref{divergencetheorem}, and the area formula \eqref{area}, recalling \eqref{div-equi} and that $U_B=\nu_{\mathbb{H}}$ on $\partial B^\lambda\setminus \R(0,0,1)$ for every $\lambda>0$: 
\begin{equation}\label{dive2}\lim_{\epsilon\to 0}\int_{B^\lambda}  \Div_g( U_B \phi_\epsilon\psi_\epsilon)  dp= \lim_{\epsilon\to 0}\int_{\partial B^\lambda} \psi_\epsilon\phi_\epsilon \nu_{\mathbb{H}} \cdot_g N_g  d\sigma=\int_{\partial B^\lambda}  |N_{\mathbb{H}}| d\sigma=\Per(B^\lambda).\end{equation}

\noindent {\bf Step 4: $B$ is outward minimizing and \eqref{outwardB} holds}. 
Let us define $\lambda_F:=  \inf\{\lambda\,|\,  F\subseteq B^\lambda\}$, with $\lambda_F=+\infty$ if such $\lambda$ does not exist. 

Recall that away from the $z$-axis we have that $\operatorname{div}_{\mathbb H}U_B(p)\geq \frac{\delta}{\rho_B(p)}$. Therefore, if $F\subseteq B^{\lambda}$, then for every $p\in F$ which is not on the $z$-axis, since $\rho_B(p)\leq \lambda$, there holds  $\operatorname{div}_{\mathbb H}U_B(p)\geq \frac{\delta}{\rho_B(p)}\geq \frac{\delta}{\lambda}$.

Let now  $F\in\mathcal{P}$ be  such that $F\supseteq B$. Then by \eqref{pansu-calibration-ineqB}, \eqref{pansu-calibration-eqB}  and recalling that $\operatorname{div}_{\mathbb H}U_B(p)\geq \frac{\delta}{\lambda_F}$ for $p\in F\setminus B$, we obtain 
\[\Per(F)\geq  \int_F \operatorname{div}_{\mathbb H}U_B(p)dp = \int_{B} \operatorname{div}_{\mathbb H}U_B(p)dp+\int_{F\setminus B} \operatorname{div}_{\mathbb H}U_B(p)dp\geq \Per(B)+\frac{\delta }{\lambda_F}|F\setminus B|.\] 

The second statement \eqref{outwardB2} follows easily  using the submodularity   \eqref{propPer}. If $F\supseteq B$, then  \eqref{outwardB2} coincides with  \eqref{outwardB}. Assume  that $F\not\supseteq B$. Let $E=F\cup B$.   Then $\lambda_E=\max\{\lambda_F,1\}$, so we  apply \eqref{outwardB} to $E$  observing that $|E\setminus B|= |F\setminus B |$ and obtain 
\begin{align*}\Per(F)+\Per(B)&\geq \Per(F\cap B)+\Per( E)\geq \Per(F\cap B)+\Per(B) + \frac{\delta}{\lambda_E}|E\setminus B|\\ &=  \Per(F\cap B)+\Per(B) + \frac{\delta}{\max\{\lambda_F, 1\}}|F\setminus B|.\end{align*}
\end{proof}

 \section{The minimizing movement scheme} 
\subsection{Definition and basic properties} \label{inizio}
We introduce the following set valued operator $T_h$ for $h>0$: given  a measurable set  $\emptyset \neq E\in \mathcal{M}$ let
\[E\rightarrow T_h(E)=\arg\min \{J_h(F)\ ,\ F \in \mathcal{M}\}\] 
where 
\begin{equation}J_h(F):= P_{\mathbb{H}} (F) + \frac{1}{h} \int_{F} d_{E}(p) dp, \label{Jh} 
\end{equation}
and  $d_E$ is  the  signed distance introduced in \eqref{sig} associated with the Carnot--Carathéodory distance $d_{CC}$  defined in \eqref{dcc}. 

\begin{remark}\upshape Note that the results presented   in  Sections \ref{inizio} and \ref{sec6}   hold for any choice of distance function presented in Section~\ref{sec:distances}. Nonetheless, from now on we will use only the Carnot--Carathéodory distance  $d_{CC}$. This choice is motivated by the analysis we will present in Section~\ref{sec:density}.
\end{remark}
Notice that in view of \eqref{propPer} the functional $J_{h}$ is submodular, meaning that
$$ J_{h}(F_{1}\cap F_{2})+ J_{h}(F_{1}\cup F_{2}) \leq J_{h}(F_{1}) + J_{h}(F_{2})$$
for any $F_{1}, F_{2} \in \mathcal{P}. $

In general the minimizers of the functional $J_h$ in \eqref{Jh} are not unique  (see Proposition \ref{propertyTh} below). With a slight abuse of notation we will denote by $T_h(E)$ a generic minimizer of $J_h$, i.e. an element $F\in T_h(E)$. 

 \begin{remark}\label{remarkdelta} \upshape
Observe that
\[
\int_{F}d_E(p)dp - \int_{E} d_E(p)dp
=
\int_{F \setminus E} d_E(p)dp - \int_{E \setminus F} d_E(p)dp= \int_{E \triangle F} dist(p, \partial E) dp.
\]
In particular 
\begin{equation}J_h^{a}(F):= P_{\mathbb{H}} (F) + \frac{1}{h} \int_{E \triangle F} dist(p, \partial E) dp= J_{h}(F)-\frac{1}{h}\int_{E} d_E(p)dp , \label{Jha} 
\end{equation}
so  $J_{h}(\cdot)$ and $J_{h}^{a}(\cdot)$ 
share the same minimizers. 
\end{remark}

We verify some  standard properties of the operator $T_h$ (cp. with \cite[Section 6]{cmp} where the euclidean case is treated).
\begin{proposition}\label{propertyTh} 
Let $E \in \mathcal{P}$ be  bounded, so that  $P_{\mathbb{H}}(E) < \infty$ and $E\subseteq B_P^R$ for some $R>0$, where $B_P^R$ is the Pansu ball of parameter $R$. 
\begin{enumerate}
\item There exists at least   one bounded  $\mathbb{H}$-Caccioppoli set  which solves the  minimization problem $\min_{F\in\mathcal{M}} J_h(F)$. Moreover  every minimizer is contained in $\overline{B_P^R}$ and  there exists a minimal $T_h^-(E)$ and a maximal solution $T_h^+(E)$ (with respect to inclusion) and $T_h^+(E) \subseteq \overline{B_{P}^{R}}$. 
  \item  If $E \subset E' \subset \mathbb{H}$  for $ E' \in \mathcal{P}$ bounded set,  then for all   $F\in T_{h}(E) $ there exists $F'\in T_{h}(E')$ such that $F\subseteq F'$.  Moreover: $T^{\pm}_{h}(E) \subset T^{\pm}_{h}(E')$.
   \item  If $F\in T_h(E)$  and $R_{\theta}E=E$ with $R_\theta$ a rotation around the $z$-axis   of angle $\theta$, then also $R_\theta F\in T_h(E)$.  Therefore, if  $E$ is rotationally symmetric with respect to the $z$-axis, i.e. $R_{\theta}E=E$
 for  every $\theta \in (0, 2\pi)$  
then $T_h^+(E)$  is rotationally symmetric with respect to the $z$-axis, that is 
$ T_h^+(E)=R_{\theta} T_h^+(E)$ for all $\theta$.   
\item  If $F\in T_h(E)$, then $\Per (F)\leq \Per(E)$.
\end{enumerate}
\end{proposition}
\begin{proof}

 \begin{enumerate}
\item  By Proposition \ref{convexset} for every $F\subseteq \mathcal{P}$ there holds $\Per(F\cap B_P^R)\leq  \Per(F)$. Moreover, since $d_E(x)>  0$ for $x\in \mathbb H\setminus B_P^R$, we obtain that  $\int_{F}d_E(p)dp\geq \int_{F\cap B_P^R} d_E(p)dp$. Therefore $J_h(F)\geq J_h(B_P^R\cap F)$, with strict inequality if $|F\setminus B_P^R|>0$.

This implies that every minimizer of $J_h$ necessarily satisfies $|F\setminus B_P^R|=0$, and so $F\subseteq\overline{B_P^R}$. Now, we are  reduced to minimize $J_h(F)$ on sets $F\subseteq \overline{B_P^R}$, $F\in\mathcal{P}$. First of all we note that $\min_{B_P^R} d_E(x)=-k\leq 0$, so for every $F\subseteq \overline{B_P^R}$, 
$
J_h(F) \geq \frac{1}{h}\int_F d_E(p)dp\geq -\frac{k|F|}{h}\geq  -\frac{k|B_P^R|}{h}.
$
 
 We consider a minimizing sequence $F_n$: in particular  for every $n$, since $d_E\geq 0$ on $\mathbb{H}\setminus E$ and $d_E\geq -k$ on $E$, it holds 
 \[ \infty> J_{h}(E)\geq J_h(F_n)\geq  \Per(F_n)+\frac{1}{h}\int_{E\cap F_n} d_E(p)dp\geq  \Per(F_n)-\frac{k}{h}|E|.\]

  Therefore  $\Per(F_n)$ are equibounded  and  
 by the compactness theorem in the space of $BV_{\mathbb{H}}$ functions   (see Theorem \ref{CompactnessStrong}), 
 we have that up to passing to a subsequence $F_n\to F$ in $L^1$ (which means, $\chi_{F_n}\to \chi_F$ in $L^1$), and $\liminf_n \Per(F_n)\geq \Per(F)$ since the perimeter functional is lower semicontinuous with respect to $L^1_{loc}$ convergence of (characteristic functions of) sets.
  
 Note that up to passing to a subsequence we have  $\chi_{F_n}\to \chi_F$  for a.e. $p$. 
 
 Then observe that by Fatou lemma,  it holds 
  \[\liminf_n\int_{\mathbb{H}\setminus E}\chi_{F_n}d_E dp\geq \int_{\mathbb{H}\setminus E} \chi_{F} d_Edp\qquad \liminf_n\int_{E}(d_E\chi_{F_n} +k)dp\geq\int_{E} (d_E\chi_{F}+k)dp\]
 Therefore we conclude that $F$ is a minimizer.

 If we consider $F_i\in T_h(E)$, for $i=1, 2$,  then $F_1\cup F_2$ and $F_1\cap F_2$ are both minimizers, by   the submodularity of $J_{h}$.

 For the last statement see \cite[Proposition~6.1]{cmp}.
 
In the rest of the proof with a slight abuse of notation we will denote as $T_h(E)$ a minimizer of $J_h$.
 \item Assume $E \subset \subset E'$, so  that $d_{E}(p)> d_{E'}(p)$ for all $p \in \mathbb{H}$. 
 
 By minimality of $T_h(E)$ we have 
\begin{equation}\label{intersection}
P_{\mathbb{H}}(T_{h}(E) \cap T_{h}(E')) + \frac{1}{h} \int_{T_{h}(E) \cap T_{h}(E')} d_E(p)dp \geq P_{\mathbb{H}}(T_{h}(E)) + \frac{1}{h} \int_{T_{h}(E)} d_E(p)dp .
\end{equation}
In a similar way, by minimality of  $T_{h}(E')$ we have 
\begin{equation}\label{union}
P_{\mathbb{H}}(T_{h}(E) \cup T_{h}(E')) + \frac{1}{h} \int_{T_{h}(E) \cup T_{h}(E')} d_{E'}(p) dp \geq P_{\mathbb{H}}(T_{h}(E')) + \frac{1}{h} \int_{T_{h}(E')} d_{E'}(p) dp.
\end{equation}
Summing the two inequalities and using the sub-additivity of the perimeter \eqref{propPer} we obtain 
\begin{equation*}
\int_{T_{h}(E) \cap T_{h}(E')} d_E(p)dp + \int_{T_{h}(E) \cup T_{h}(E')} d_{E'}(p) dp \geq \int_{T_{h}(E)} d_E(p)dp + \int_{T_{h}(E')} d_{E'}(p) dp
\end{equation*}
from which we deduce that
\begin{equation*}
\int_{T_{h}(E) \setminus T_{h}(E')} \big(d_{E'}(p) - d_{E}(p) \big) dp \geq 0
\end{equation*}
but, since $d_{E'}(p) - d_{E}(p)<0$ then $T_{h}(E) \setminus T_{h}(E') = \emptyset$. Hence   $T_{h}(E) \subset T_{h}(E')$.


The statement for $E \subset E'$ is achieved by a perturbation argument as explained in \cite[Lemma~6.2]{cmp}.


\item By Remark \ref{inv}, we know that 
$$ P_{\mathbb{H}}(F)= P_{\mathbb{H}}(R_{\theta}F). $$
Moreover,   we have that 
\begin{align*}
\int_{R_{\theta} F} d_E(p)dp 
&=- \int_{E}\chi_{R_{\theta} F}(p) \text{dist}(p, \mathbb{H} \setminus E) dp+ \int_{\mathbb{H} \setminus E}\chi_{R_{\theta} F}(p) \text{dist}(p,E) dp \\
&=- \int_{E}\chi_{ F}(q) \text{dist}(R_{\theta}q, R_{\theta} ( \mathbb{H} \setminus E)) dq + \int_{\mathbb{H}^{1} \setminus E}\chi_{ F}(q) \text{dist}(R_{\theta}q, R_{\theta}E) dq  \\
&=\int_{ F} d_{E}(q) dq
\end{align*}
where we have used that
$ \text{dist}(R_{\theta}p, R_{\theta}q)= \text{dist}(p,q).$ 
 
 \item  Using minimality of  $T_h(E)$ and $E$ as a competitor we obtain 
\[
\Per (T_h(E)) + \frac{1}{h} \int_{T_h(E)} d_E(p)dp \leq P_{\mathbb{H}} (E) + \frac{1}{h} \int_E d_E(p)dp.
\]
Observe that 
\[\frac{1}{h} \int_E d_E(p)dp- \frac{1}{h} \int_{T_h(E)} d_E(p)dp= \int_{E\setminus T_h(E)} d_{E}(p)dp - \int_{T_h(E)\setminus E}   d_{E }(p) dp \leq 0.
\]
from which we conclude. 
 
 \end{enumerate} 
  \end{proof}

  We provide a description of the scheme for outward minimizing sets. 
  \begin{proposition}\label{outward}
 Let  $E \in \mathcal{P}$  be a bounded outward minimizing in $\mathbb{H}$ with $|\partial E|=0$. Then 
\begin{enumerate}
  \item    $T_h(E)\subseteq E$, that is every minimizer $F\in T_h(E)$ satisfies $F\subseteq E$, 
\item  every minimizer in $T_h(E)$ is   outward minimizing in $\mathbb{H}$. 
\end{enumerate} 
  \end{proposition} \begin{proof} Note that if $|\partial E|=0$, then $d_E(p)>0$ for a.e. $p\in \mathbb H\setminus E$.
  \begin{enumerate}

\item Let  $F$ be a minimizer of $J_h$ in \eqref{Jh}. If $E$ is an outward minimizing set,  then by the characterization in Lemma~\ref{lemma1.1} we have that $\Per(F\cap E)\leq \Per(F)$. Moreover
\[\int_{F\cap E} d_E(p)dp\leq \int_{F} d_E(p)dp\] with strict inequality if  $F\setminus  E$ has positive measure. Therefore if $|F\setminus E|>0$,  it holds that
$J_h(F\cap E)<J_h(F)$, which implies that $F$ cannot be a minimizer of $J_h$. This implies that necessarily every minimizer $T_h(E)$ is contained in $E$. 

\item Let   $G\in\mathcal{P}$ with $F\subset G$, where $F \in T_{h}(E)$.
  Then, using $G\cap E$ as a competitor we get, recalling that $F\subseteq E\cap G$ and $d_E\leq 0$ in $E$, 
   \begin{eqnarray*}  P_{\mathbb{H}} (F) &\leq&  P_{\mathbb{H}} ( E\cap G) +  \frac{1}{h} \int_{E\cap G} d_E(p)dp- \frac{1}{h} \int_{F} d_E(p)dp  
\\ &=&
 P_{\mathbb{H}} (E\cap G) +  \frac{1}{h} \int_{E\cap G \setminus F} 
 d_E(p)dp\leq P_{\mathbb{H}} (E\cap G)\leq P_{\mathbb{H}} (G) \end{eqnarray*}  where the last inequality is a consequence of the fact that 
 $E$ is outward minimizing in $\mathbb{H}$ and of the characterization proved in Lemma \ref{lemma1.1}. 

  \end{enumerate} 
  \end{proof}
  

  We present an  estimate on the one time step of the scheme for   sets $B$ in the class $\mathcal{B}$ introduced in Definition \ref{definitionb} and which are also horizontal strictly mean convex.
The following Proposition~\ref{speedpansu} applies in particular to Pansu balls. Observe that $T_{h}(B)$ might  no longer belong to $\mathcal{B}$.
For $\lambda >0$ recall that  $ B^{\lambda}:= \delta_{\lambda}(B). $
 \begin{proposition}\label{speedpansu} Let $B\in \mathcal{B}$ be horizontal  strictly mean convex with constant $\delta>0$, and let  \[  \lambda^0=\inf\{\lambda>0 \text{ s.t. } |  T_h(B)\setminus B^{\lambda }|=0\}\} .\] 
   Then  there exist $h_0=h_0(diam(B), \delta)$, $C=C(diam(B), \delta)$, such that   for all $0<h\leq h_0$ there holds 
    \[  \lambda^0\leq \sqrt{1-Ch^{2} } .\]
   \end{proposition}
   \begin{proof}
First of all we observe that  by Proposition \ref{convexset} the set $B$ is outward minimizing, therefore by Proposition \ref{outward} we have that $T_h(B)\subseteq B$. 
 Let  $\lambda\in (0,1 )$ and  $  B^{\lambda} \subseteq B$: note that  $B^{\lambda}$ is still open, bounded, $C^{2}$, axially symmetric and horizontal mean convex by Remark~\ref{remarkhomoH} with $H_{\mathbb{H}}(p, B^{\lambda}) \geq \frac{\delta}{\lambda}$. 
  Since $T_h(B)\subseteq B$ obviously $\lambda^{0} \leq 1$.
    Using Proposition \ref{convexset} applied to $B^{\lambda}$ and to $T_h(B)$, and the fact that $T_h(B)\subseteq B$, we obtain that
 \begin{equation}\label{delta1bis}   \Per( T_h(B)\cap  B^{\lambda} )\leq \Per( T_h(B))-\delta  |T_h(B)\setminus  B^{\lambda}|.\end{equation} 
 By  minimality of $T_h(B)$, using as a competitor $ T_h(B)\cap  B^{\lambda}$ we obtain 
 \begin{equation}\label{delta2bis}   \Per( T_h(B)  )\leq \Per( T_h(B)\cap  B^{\lambda})-\frac{1}{h} \int_{T_h(B)\setminus  B^{\lambda}}d_{B}(p)dp.\end{equation}
Then by \eqref{delta1bis}, \eqref{delta2bis} we infer that
\begin{equation}\label{delta3bis}
\delta |T_h(B)\setminus  B^{\lambda}|\leq -\frac{1}{h} \int_{T_h(B)\setminus  B^{\lambda}}d_{B}(p)dp
=\frac{1}{h} \int_{T_h(B)\setminus  B^{\lambda}}\text{dist}(p,\partial B) dp \end{equation} 
where in the last equality we used   the fact that $T_h(B)\subseteq B$. \\
Note that  if $p\in T_h(B)\setminus B^{\lambda}$ then $\delta_{1/\lambda} p\not \in B$, and so, recalling \eqref{dcclambda} we have \[
\text{dist}(p,\partial B) \leq d_{CC}(p, \delta_{1/\lambda} p)\leq  C\, diam(B) \max\left(|\lambda^{-1}-1|, \sqrt{|\lambda^{-2}-1|} \right).\]
We substitute in \eqref{delta3bis} and we get, for every $\lambda\in (0,1)$ such that  $ |T_h(B)\setminus B^{\lambda}|\neq 0$, 
\begin{equation}\label{conclusionebis}\delta \leq \frac{C diam(B)}{h} \max\left(|\lambda^{-1}-1|, \sqrt{|\lambda^{-2}-1|} \right)=\frac{C diam(B)}{h\lambda}\max\left((1-\lambda), \sqrt{(1-\lambda^2)}\right). \end{equation}  
Note that if $\lambda < \lambda^{0}$ then $ |T_h(B)\setminus B^{\lambda}|\neq 0$ and \eqref{conclusionebis} holds. 
We have two possibilities:
\begin{itemize}\item[(i)$\lambda^{0}\leq \frac{1}{2}$:]
In this case  \[T_h(B)\subseteq B^{1/2}\] and we conclude. 
\item[(ii)$\lambda^{0}>\frac{1}{2}$:] For $\frac{1}{2}<\lambda<\lambda^0$, \eqref{conclusionebis} reads  
\begin{equation}\label{conclusione2bis}  \sqrt{\frac{3}{4}} \geq \sqrt{(1-\lambda^2)}=  \max\left((1-\lambda), \sqrt{(1-\lambda^2)}\right)\geq \frac{h \delta}{2C diam (B)}. \end{equation}
In light of the previous inequality we set $h_{0} < 2\sqrt{\frac{3}{4}}C \frac{ diam(B)}{\delta} $. 
 Furthermore we conclude by letting $\lambda \to \lambda^{0}$ in \eqref{conclusione2bis}  that $   \sqrt{(1-(\lambda^0)^2)}\geq \frac{h \delta }{2C diam(B)}$, that is $(\lambda^0)^{2}\leq 1-\frac{h^2 \delta^{2}}{4 C^2 diam(B)^{2}}$.  
  \end{itemize} 
  \end{proof}

\subsection{The scheme for unbounded sets}\label{sec6}
 Here we extend the time discrete scheme to the case of unbounded sets $E$ with compact boundaries.
 More generally  we assume that the set $E$ in \eqref{Jh} is such that $P_{\mathbb{H}}(E) <\infty$ and   $E^{c}$ is bounded. Recalling \eqref{eq30}, the idea is to define the evolution as the complement of the evolution of the complementary set.
More precisely, using Proposition~\ref{propertyTh}, we set
 \begin{align*}
 T^{\pm}_{h}(E):=\mathbb{H} \setminus T_{h}^{\mp}(E^{c}).
 \end{align*}
 In particular note that if  $E_{1}$, $E_{2} \in \mathcal{M}$ are such that $P_{\mathbb{H}}(E_{i}) <\infty$, and $E_{i}^{c}$ are bounded, then it follows from Proposition~\ref{propertyTh}(2.) that
 $$E_{1} \subset E_{2} \Longrightarrow  T^{\pm}_{h}(E_{1}) \subset T^{\pm}_{h}(E_{2}).$$
 \begin{lemma}
 Let $E_{1}$, $E_{2} \in \mathcal{M}$ be such  that $P_{\mathbb{H}}(E_{i})<\infty$, $i=1,2$. If $E_{1} \subset E_{2}$ and $ E_{2}^{c} \cup E_{1}  \subset B^{R}_{P}$ for a sufficiently large $R>0$, then
 $T^{\pm}_{h}(E_{1}) \subset T^{\pm}_{h}(E_{2})$.
 \end{lemma}
 \begin{proof}
 Since Pansu balls are outward minimizing  sets we have by Proposition~\ref{outward} (and Proposition~\ref{propertyTh}) that
 $$ T_{h}^{\pm}(E_{2}^{c}) \subset T_{h}^{\pm}(B^{R}_{P}) \subset B^{R}_{P}. $$ 
  The proof now follows the arguments given in  \cite[Lemma~6.10]{cmp}, which we report here for completeness. Since $T^{-}_{h}(E_{2}^{c})$ is the minimal solution of
 \begin{align*}
 \min \left \{ \Per( F) +\frac{1}{h} \int_{F} d_{E_{2}^{c}} dp \, : \, F \in \mathcal{M} \right \},
 \end{align*}
 we have (using the change of variable $\tilde{F}:= \mathbb{H} \setminus F=F^{c}$, the fact that $\Per(\mathbb{H} \setminus F)=\Per( F)$, and $d_{E_{2}^{c}}=-d_{E_{2}}$) that $T^{+}_{h}(E_{2})= \mathbb{H} \setminus T_{h}^{-}(E_{2}^{c})$ is the maximal solution of
\[\min \left \{ \Per( \tilde{F}) +\frac{1}{h}\int_{B^{R}_{P}}d_{E_{2}} dp-\frac{1}{h} \int_{\tilde{F}^{c}} d_{E_{2}} dp \, : \, \tilde{F} \in \mathcal{M} \right \} - \frac{1}{h}\int_{B^{R}_{P}}d_{E_{2}} dp.\]
 Note that if $ \tilde{F}^{c} \subset B^{P}_{R}$, then $\int_{B^{R}_{P}}d_{E_{2}} dp- \int_{\tilde{F}^{c}} d_{E_{2}} dp= \int_{\tilde{F}}d_{E_{2}} \chi_{B^{R}_{P}} dp$. Therefore, since $T^{+}_{h}(E_{2})^c\subset B^{P}_{R}$, we conclude that $T^{+}_{h}(E_{2})$ is  also the maximal solution of
 \begin{align}\label{eq69}
 \min \left \{ \Per( \tilde{F}) +\frac{1}{h}\int_{\tilde{F}}d_{E_{2}} \chi_{B^{R}_{P}} dp \, : \, \tilde{F} \in \mathcal{M} , \,\tilde{F}^{c} \subset B^{P}_{R}\right \}.
 \end{align}
Analogously one can show that $T^{-}_{h}(E_{2})$ is the minimal solution of \eqref{eq69}.\\
 Since $d_{E_{2}} \chi_{B^{R}_{P}}  \leq d_{E_{1}}$ and since $T^{\pm}_{h}(E_{1}) \cup T^{\pm}_{h}(E_{2})$, $T^{\pm}_{h}(E_{1}) \cap T^{\pm}_{h}(E_{2})$  are admissible sets in \eqref{eq69}, one can argue as in Proposition~\ref{propertyTh} (2.) to infer that $T^{\pm}_{h}(E_{1}) \subset T^{\pm}_{h}(E_{2}) $.
 \end{proof}

\subsection{Density estimates}\label{sec:density}
In this section, we provide density estimates for minimizing solutions to \eqref{Jh}.
 The arguments are quite classical (we refer to \cite{lu}, \cite{Mugnai}, \cite{ChamDeG} and reference therein) and we provide sketch of the main proofs with adaptations to our setting. 

We start  recalling the following important result, which provides a relation between the $\mathbb{H}$-perimeter of a set $F$ on the boundary of a ball and the derivative of its volume.
\begin{lemma}{(\cite[Lemma~1.19]{lr})}\label{lemma5.1}
Let $F$ be a set with finite $\mathbb{H}$-perimeter,    $p \in \mathbb{H}$ fixed, and  define 
\begin{equation}\label{eq53}
\begin{cases}  P_{\mathbb{H}}(F \cap B_{CC}(p, r), \partial B_{CC}(p, r)) := P_{\mathbb{H}}(F \cap B_{CC}(p, r)) - P_{\mathbb{H}}(F, B_{CC}(p, r))\\ 
 P_{\mathbb{H}}(F \setminus B_{CC}(p, r), \partial B_{CC}(p, r)):= P_{\mathbb{H}}(F \setminus B_{CC}(p, r)) - P_{\mathbb{H}}(F, \mathbb{H} \setminus \overline{B_{CC}}(p,r)). \end{cases}\end{equation} 
For $r>0$ let   $m(r):=|F \cap B_{CC}(p, r)|$. Then  for  a.e.  $r >0$,  there holds \begin{align*}
P_{\mathbb{H}}(F \cap B_{CC}(p, r), \partial B_{CC}(p, r)) \leq m'(r),\\
P_{\mathbb{H}}(F \setminus B_{CC}(p, r), \partial B_{CC}(p, r)) \leq m'(r).
\end{align*}
 
\end{lemma}
 
We recall a technical  Lemma (see e.g. \cite[Lemma~4.4]{Mugnai}).  We provide a short sketch of the proof. 
\begin{lemma}\label{lemma9.4}
 Let $F$ be a set of finite $\mathbb{H}$-perimeter.
\begin{itemize} \item [(i)]
Assume that  $|B_{CC}(\hat{p},r) \setminus F| >0$ for all $r >0$ and that there exist $R>0$, $\mu\geq 0$ such that 
\begin{align}\label{eq84}
P_{\mathbb{H}}(F) \leq P_{\mathbb{H}}(F \cup B_{CC}(\hat{p},r)) +\mu  |B_{CC}(\hat{p},r) \setminus F| \quad \forall\quad  0<r<R.
\end{align}
Then, there exist   $\sigma>0, C>0$ depending only on the isoperimetric constant \eqref{isoineq} such that 
\begin{align}\label{eq85}
 |B_{CC}(\hat{p},r) \setminus F| \geq \sigma r^{4}  \quad \forall \quad 0<r< \min \left\{ R, \frac{C}{\mu} \right \}. 
\end{align}
\item[(ii)] Assume that  $|B_{CC}(\hat{p},r) \cap F| >0$ for all $r >0$ and that there exist $R>0$, $\mu\geq 0$ such that  
\begin{align}\label{eq842}
P_{\mathbb{H}}(F) \leq P_{\mathbb{H}}(F \setminus B_{CC}(\hat{p},r)) +\mu  |B_{CC}(\hat{p},r) \cap F| \quad \forall\quad  0<r<R.
\end{align}
Then,  there exist   $\sigma>0, C>0$ depending only on the isoperimetric constant \eqref{isoineq} such that 
\begin{align}\label{eq852}
 |B_{CC}(\hat{p},r) \cap F| \geq \sigma r^{4}  \quad \forall \quad 0<r< \min \left\{ R, \frac{C}{\mu} \right \}.
\end{align}
\end{itemize} For $\mu=0$ in \eqref{eq842}, \eqref{eq852}, we set $\min \left\{ R, \frac{C}{\mu} \right \}=R$. \end{lemma} 
\begin{proof}
(i) 
Recalling  \eqref{eq30}  and \eqref{eq53} for $F^c$, since  $\mathbb{H}\setminus (F^c\cap B_{CC}(\hat{p},s))=F\cup B_{CC}(\hat{p},s)$ and $F^c\cap B_{CC}(\hat{p},s) = B_{CC}(\hat{p},s) \setminus F$, we   observe that
\begin{align}\label{yesy2}
&P_{\mathbb{H}}(F \cup B_{CC}(\hat{p},s))= P(F, \mathbb{H}\setminus \overline{B_{CC}(\hat{p},s)}) 
 + P_{\mathbb{H}}(  B_{CC}(\hat{p},s)\setminus F, \partial B_{CC}(\hat{p},s)) \notag\\
 \leq& P_{\mathbb{H}}(F, \mathbb{H})
- P_{\mathbb{H}}(F, B_{CC}(\hat{p},s)) +P_{\mathbb{H}}(  B_{CC}(\hat{p},s) \setminus F, \partial  B_{CC}(\hat{p},s))\notag 
 \\
=&P_{\mathbb{H}}(F)
- P_{\mathbb{H}}(F^{c}, B_{CC}(\hat{p},s)) +P_{\mathbb{H}}( B_{CC}(\hat{p},s) \setminus F, \partial  B_{CC}(\hat{p},s))\notag \\
=& P_{\mathbb{H}}(F ) -  P_{\mathbb{H}}( F^c\cap B_{CC}(\hat{p},s)) +P_{\mathbb{H}}( F^c\cap B_{CC}(\hat{p},s), \partial  B_{CC}(\hat{p},s))+P_{\mathbb{H}}( B_{CC}(\hat{p},s) \setminus F, \partial  B_{CC}(\hat{p},s))\notag 
\\ =&P_{\mathbb{H}}(F ) -  P_{\mathbb{H}}( B_{CC}(\hat{p},s) \setminus F) + 2 P_{\mathbb{H}}(B_{CC}(\hat{p},s) \setminus F, \partial  B_{CC}(\hat{p},s)). 
\end{align}  
From \eqref{yesy2} and \eqref{eq84} we obtain for all $s<R$, 
\[
P_{\mathbb{H}}( B_{CC}(\hat{p},s) \setminus F) \leq 2 P_{\mathbb{H}}( B_{CC}(\hat{p},s) \setminus F, \partial  B_{CC}(\hat{p},s)) + \mu |B_{CC}(\hat{p},s) \setminus F|.
\]
We set $|B_{CC}(\hat{p},s) \setminus F|=|B_{CC}(\hat{p},s) \cap F^{c}|=:m(s)$. So, by Lemma~\ref{lemma5.1},  and the isoperimetric inequality \eqref{isoineq} we get for  a.e. $s <R$: 
\[ C m(s)^{\frac{3}{4}} \leq 2 m'(s) +\mu m(s).\]   Since 
$0<m(s) \leq  |B_{CC}(0,1)|s^{4}$ we get $m(s) \leq |B_{CC}(0,1)|^{\frac{1}{4}} s  \,m(s)^{1-\frac{1}{4}}$: then substituting in the previous inequality we conclude   
$$\frac{C}{2}m(s)^{\frac{3}{4}} \leq \left(C-s\mu |B_{CC}(0,1)|^{\frac{1}{4}}\right) m(s)^{\frac{3}{4}} \leq 2 m'(s)$$
for  $s \leq \min \{ R,\frac{C}{2 |B_{CC}(0,1)|^{\frac{1}{4}}} \frac{1}{\mu} \} $, when $\mu>0$, for $s\leq R$ when $\mu=0$. 
By direct integration, we get the claim.

(ii) 
Recalling \eqref{eq53}, we have that 
\begin{align}\label{yesy}
&P_{\mathbb{H}}(F\setminus B_{CC}(\hat{p},s)) =P_{\mathbb{H}}(F,\mathbb{H} \setminus \overline{B_{CC}(\hat{p},s)}) +P_{\mathbb{H}}(F \setminus B_{CC}(\hat{p},s), \partial  B_{CC}(\hat{p},s))
\notag\\
& \leq P_{\mathbb{H}}(F,\mathbb{H}) - P_{\mathbb{H}}(F, B_{CC}(\hat{p},s)) +P_{\mathbb{H}}(F \setminus B_{CC}(\hat{p},s), \partial  B_{CC}(\hat{p},s))\notag \\
 &= P_{\mathbb{H}}(F) -P_{\mathbb{H}}(F \cap B_{CC}(\hat{p},s)) +P_{\mathbb{H}}(F \cap B_{CC}(\hat{p},s), \partial  B_{CC}(\hat{p},s)) \notag\\& \quad 
+P_{\mathbb{H}}(F \setminus B_{CC}(\hat{p},s), \partial  B_{CC}(\hat{p},s)).
\end{align}
The claim follows in a similar way to part (i),
 using now \eqref{yesy} and \eqref{eq842}.  
 
\end{proof}

  The next  result  is an adaptation to our setting of \cite[Lemma~3.5]{ChamDeG}, and it is an estimate  independent of the time step $h$. 
\begin{lemma}\label{lemma8.1b}
Let $E \in \mathcal{P}$ be bounded. Let $F\in T_{h}(E)$.
Then there exists a constant $\sigma >0$ (depending only on the isoperimetric constant in \eqref{isoineq}) with the following property: 
\begin{itemize}\item [(i)] if 
$|F \cap B_{CC}(\hat{p},s)|>0$ for all $s>0$ and  $E \cap B_{CC}(\hat{p},r_0)=\emptyset$ for some $r_0$
then
\begin{align}\label{eq-74b}
|F \cap B_{CC}(\hat{p},s)| \geq \sigma s^{4} \qquad \forall s\leq r_0.
\end{align}
\item [(ii)] if 
$| B_{CC}(\hat{p},s) \setminus F|>0$ for all $s>0$ and $B_{CC}(\hat{p},r_0) \subset E$ for some $r_0$, 
then
\begin{align}\label{eq-75b}
| B_{CC}(\hat{p},s) \setminus F|\geq \sigma s^{4}\qquad \forall s\leq r_0.
\end{align}\end{itemize}
\end{lemma}
\begin{proof}
(i)   From the minimality of $F$ we obtain for any $s >0$ that
\begin{align*}
P_{\mathbb{H}}(F)  +\frac{1}{h} \int_{F} d_{E}(p)\, dp \leq P_{\mathbb{H}}(F\setminus B_{CC}(\hat{p},s))  +\frac{1}{h} \int_{F\setminus B_{CC}(\hat{p},s)} d_{E}(p) \,dp.
\end{align*}
Using that $d_{E}$ is positive on $F \cap B_{CC}(\hat{p},s) \subset E^{c}$   we infer that
$$P_{\mathbb{H}}(F) \leq  P_{\mathbb{H}}(F\setminus B_{CC}(\hat{p},s)).$$
The claim now follows from Lemma~\ref{lemma9.4}(ii) with $\mu=0$.

To prove (ii) one argues in a similar way: first use the minimality of $F$ and test against $F \cup B_{CC}(\hat{p},s) $. Since $d_{E}$ is negative on $B_{CC}(\hat{p},s)  \subset E$ (for $s \in (0,r_0)$) we obtain
$$ P_{\mathbb{H}}(F) \leq P_{\mathbb{H}}(F \cup B_{CC}(\hat{p},s)).$$  Application of Lemma~\ref{lemma9.4}(i) with $\mu=0$ yields the claim.
\end{proof}
  
We provide some uniform bounds between the boundary of a minimizing set and the boundary of the reference set. 
\begin{proposition}($L^{\infty}-estimate$)\label{inftyesti} 
There exists a constant $\gamma>0$ depending only on the isoperimetric constant  such that for any  
 bounded  $E \in \mathcal{P}$ and $F \in T_{h}(E)$ there holds
 \begin{align}\label{eq88}
 \sup_{E \triangle F} |d_{E}| \leq \gamma \sqrt{h}.
 \end{align}
  \end{proposition}
  \begin{proof} The proof follows  \cite[Proposition 3.2]{Mugnai}. We set  $\gamma=\max\left\{\frac{4P_{\mathbb{H}}(B_{CC}(0, 1)) }{\sigma} , 2\right\}$,  where $\sigma$ is the constant depending only on the isoperimetric constant appearing in Lemma \ref{lemma9.4}.  
  
  We argue by contradiction.
  
   Let $\hat{p} \in E \setminus  F$, and assume that  
 \begin{align}\label{eq89}
 d_{E}(\hat{p}) < -\gamma \sqrt{h}.
 \end{align} 
This implies that $B_{CC}(\hat{p}, r) \subset E$ for all $r \leq \frac{\gamma}{2}\sqrt{h}$.
  Using the minimality of $F$ 
  we obtain
  \begin{align*}
  P_{\mathbb{H}}(F) \leq P_{\mathbb{H}}(F \cup B_{CC}(\hat{p}, r)) +\frac{1}{h} \int_{B_{CC}(\hat{p}, r) \setminus F} d_{E}(p) \, dp
  \end{align*}
  From \eqref{eq89} and $r \leq \frac{\gamma}{2}\sqrt{h}$ we infer that
 $$ \frac{1}{h} \int_{B_{CC}(\hat{p}, r) \setminus F} d_{E}(p) \, dp \leq \frac{1}{h} (r-\gamma  \sqrt{h})| B_{CC}(\hat{p}, r) \setminus F|\leq - \frac{\gamma}{2\sqrt{h}}| B_{CC}(\hat{p}, r) \setminus F|$$
 so that we obtain
 \begin{align}\label{eq90}
 P_{\mathbb{H}}(F) \leq P_{\mathbb{H}}(F \cup B_{CC}(\hat{p}, r)) -\frac{\gamma}{2 \sqrt{h}}| B_{CC}(\hat{p}, r) \setminus F|.
 \end{align}
 Since $\hat p\not\in F$, then $|B_{CC}(\hat{p}, r) \setminus F|>0$ for all $r$ and  we may apply   Lemma~\ref{lemma9.4} (i) with $\mu=0$ to get  \begin{align*}
 | B_{CC}(\hat{p}, r) \setminus F| \geq \sigma r^{4} \qquad \text{ for  } 0<r<\frac{\gamma\sqrt{h}}{2}.
 \end{align*}
 On the other hand, from \eqref{eq90} and \eqref{propPer} we have
 \begin{align*}
 \frac{\gamma}{2 \sqrt{h}}| B_{CC}(\hat{p}, r) \setminus F| &\leq P_{\mathbb{H}}(F \cup B_{CC}(\hat{p}, r)) -P_{\mathbb{H}}(F)  \\
 &\leq P_{\mathbb{H}}(B_{CC}(\hat{p}, r)) - P_{\mathbb{H}}(F \cap B_{CC}(\hat{p}, r)) \\
 &\leq P_{\mathbb{H}}(B_{CC}(\hat{p}, r))= P_{\mathbb{H}}(B_{CC}(0, 1))r^{3}.
 \end{align*} 
 Combining the last two inequalities and recalling the choice of $\gamma$, we get  for  $0<r<\frac{\gamma\sqrt{h}}{2}$
 $$ \sigma r^{4} \leq | B_{CC}(\hat{p}, r) \setminus F| \leq  P_{\mathbb{H}}(B_{CC}(0, 1)) \frac{2\sqrt{h}}{\gamma}r^{3} \leq \frac{\sigma\sqrt{h}}{2 }r^{3}   
$$ 
which gives a contradiction as soon as $\frac{ \sqrt{h}}{2}\leq \frac{\gamma\sqrt{h}}{4}<r < \frac{\gamma\sqrt{h}}{2}$. 

The case $\hat{p} \in F\setminus E$ is treated in a similar way.  Assume by contradiction that $|d_{E}(\hat{p}) |=d_{E}(\hat{p}) >\gamma \sqrt{h}$, so that  $B_{CC}(\hat{p}, r) \subset \mathbb{H}\setminus E$ for any $r \leq \frac{\gamma}{2}\sqrt{h}$.
Using the minimality of $F$ and testing against $F \setminus B_{CC}(\hat{p}, r)$ we obtain for  $r \leq \frac{\gamma}{2}\sqrt{h}$:
\begin{align*}
  P_{\mathbb{H}}(F) \leq P_{\mathbb{H}}(F \setminus B_{CC}(\hat{p}, r)) -\frac{1}{h} \int_{F \cap B_{CC}(\hat{p}, r) } d_{E}(p) \, dp \leq P_{\mathbb{H}}(F \setminus B_{CC}(\hat{p}, r)) -\frac{\gamma}{2 \sqrt{h}}|F \cap B_{CC}(\hat{p}, r)|.
  \end{align*}
  Using Lemma~\ref{lemma9.4}(ii) with $\mu=0$  we have 
  \begin{align*}
  |F \cap B_{CC}(\hat{p}, r)| \geq  \sigma r^{4} \qquad \text{for }  0<r<\frac{\gamma}{2}\sqrt{h}.
  \end{align*}
  On the other hand,  as above  by the previous inequality and \eqref{propPer},  \begin{align*}
  \frac{\gamma}{2 \sqrt{h}}|F \cap B_{CC}(\hat{p}, r)| &\leq P_{\mathbb{H}}(F \setminus B_{CC}(\hat{p}, r)) -P_{\mathbb{H}}(F)=P_{\mathbb{H}}(F^c\cup   B_{CC}(\hat{p}, r)) -P_{\mathbb{H}}(F^c) \\
   &\leq P_{\mathbb{H}}(B_{CC}(\hat{p}, r)) - P_{\mathbb{H}}(F^c \cap B_{CC}(\hat{p}, r))\leq  P_{\mathbb{H}}(B_{CC}(0, 1))r^{3}.  \end{align*}
 As before, we get a contradiction as soon as $\frac{ \sqrt{h}}{2}\leq \frac{\gamma\sqrt{h}}{4}<r < \frac{\gamma\sqrt{h}}{2}$.
  \end{proof}
  Finally we provide  density estimates, using Lemma \ref{lemma8.1b} and the $L^\infty$ estimate \eqref{eq88}. 
  \begin{proposition}\label{density} 
  Let $E \in \mathcal{P}$  be bounded  and let $F \in T_{h}(E)$. 
  
  There exists positive constants $\sigma$, $\theta$, $\Theta$, depending only on the isoperimetric constant \eqref{isoineq} such that the following holds:
  for any  $\hat{p}\in \partial^{e} F$ there holds
  \begin{align}\label{dea1}
  \min \{ |B_{CC}(\hat{p},r) \setminus F|, |B_{CC}(\hat{p},r) \cap F| \} \geq \sigma r^{4},\\ \label{dea2}
  \theta r^{3} \leq P_{\mathbb{H}}(F, B_{CC}(\hat{p},r)) \leq \Theta  r^{3}
  \end{align} for all   $r \in (0, C \sqrt{h})$, where $C$  only depends on the isoperimetric constant.
  \end{proposition}
  
  \begin{remark} \label{rem-open-close}\upshape  Note that by the density estimate,  for every $F\in T_h(E)$,  
  \[F(1)=\{p\in \mathbb{H}, \exists r>0, |B_{CC}(p, r)\setminus F|=0\}\qquad  F(0)=\{p\in \mathbb{H}, \exists r>0, |B_{CC}(p, r)\cap F|=0\},\] so they  are open set, moreover there holds that $|\partial^e F|=0$ and $\partial^e F=\partial F(1)=\partial F(0)$, see \cite[Theorem 3.4]{lr}. Therefore   for every $F\in T_h(E)$ we have an open representative which is $F(1)$ and a closed representative which is $\overline{ F(1)}$. 
  From now on we choose  the maximal  solution of \eqref{Jh} to be closed and the minimal to be open. 
\end{remark}
  \begin{proof}
   Since $\hat{p} \in \partial^{e}F = \mathbb{H} \setminus (F(1) \cup F(0))$ we have that 
   $$  \min \{ |B_{CC}(\hat{p},r) \setminus F|, |B_{CC}(\hat{p},r) \cap F| \} >0 \qquad  \forall r>0,$$
   so that Lemma~\ref{lemma9.4} can be applied.
   Moreover, by \eqref{eq88} we have that $|d_E(\hat p)| \leq \gamma\sqrt{h}$ (indeed if for some $\hat{p} \in \partial^{e}F$, by contradiction we have $|d_E(\hat p)| > \gamma\sqrt{h}$, then in every neighborhood of $\hat p$ we would have a subset of $E\Delta F$ of positive measure where  \eqref{eq88} does not hold). 
   
  Let us fix $r\leq \gamma\sqrt{h}$. 
Now we use  the minimality of $F$, testing against the competitor $F \cup B_{CC}(\hat{p},r) $ and  $ F\setminus B_{CC}(\hat{p},r) $  to obtain    
\[  \begin{cases}P_{\mathbb{H}}(F) 
  & \leq  P_{\mathbb{H}}(F\cup B_{CC}(\hat{p}, r)) +\frac{1}{h}\int_{  B_{CC}(\hat{p}, r)\setminus F } d_E dp\leq  P_{\mathbb{H}}(F\cup B_{CC}(\hat{p}, r)) +\frac{2\gamma}{\sqrt{h}} |B_{CC}(\hat{p}, r)\setminus F| \\
 P_{\mathbb{H}}(F) &
  \leq  P_{\mathbb{H}}(F\setminus B_{CC}(\hat{p}, r)) -\frac{1}{h}\int_{  B_{CC}(\hat{p}, r)\cap F } d_E dp\leq  P_{\mathbb{H}}(F\setminus  B_{CC}(\hat{p}, r)) +\frac{2\gamma}{\sqrt{h}} |B_{CC}(\hat{p}, r)\cap F|.\end{cases}  \]  
   
  Therefore by Lemma~\ref{lemma9.4} with $\mu=\frac{2\gamma}{\sqrt{h} }$, we conclude \eqref{dea1} for every $r<\frac{C}{2\gamma}\sqrt{h}$.
  
   The first inequality in \eqref{dea2} follows from Theorem~\ref{teoBallIso}.
 
 To show the last inequality we proceed as follows.
 From \eqref{eq53} we deduce
 \begin{align*}
 P_{\mathbb{H}}(F \cup B_{CC}(\hat{p}, r)) &= P_{\mathbb{H}}(F, \mathbb{H} \setminus \overline{B_{CC}(\hat{p}, r)}) +  P_{\mathbb{H}}(F \cup B_{CC}(\hat{p}, r), \partial B_{CC}(\hat{p}, r))\\
 & \leq  P_{\mathbb{H}}(F, \mathbb{H}) -P_{\mathbb{H}}(F, B_{CC}(\hat{p}, r))  +  P_{\mathbb{H}}(F \cup B_{CC}(\hat{p}, r), \partial B_{CC}(\hat{p}, r)).
 \end{align*} 
 and from the inequality we conclude that
 \begin{align*}
 P_{\mathbb{H}}(F, B_{CC}(\hat{p}, r)) \leq  P_{\mathbb{H}}(F \cup B_{CC}(\hat{p}, r), \partial B_{CC}(\hat{p}, r)) +\frac{2\gamma}{\sqrt{h}}|B_{CC}(\hat{p}, r) \setminus F|.
 \end{align*}
 On the other hand by \eqref{eq30} and Lemma~\ref{lemma5.1} we deduce that for a.e. $r>0$:
 \begin{align*}&
 P_{\mathbb{H}}(F \cup B_{CC}(\hat{p}, r), \partial B_{CC}(\hat{p}, r)) =P_{\mathbb{H}}(F^{c} \setminus B_{CC}(\hat{p}, r), \partial B_{CC}(\hat{p}, r))\\
& \leq \frac{d}{dr} \left( |F^{c} \cap B_{CC}(\hat{p}, r) |\right) =\frac{d}{dr} \Big( | B_{CC}(\hat{p}, r) | - | (F \cap B_{CC}(\hat{p}, r)) |\Big)\\
& =4r^{3}|B_{CC}(0, 1)| -\frac{d}{dr} \left( | (F \cap B_{CC}(\hat{p}, r)) |\right) \leq 4r^{3}|B_{CC}(0, 1)|
 \end{align*}
 since the map $r \to | (F \cap B_{CC}(\hat{p}, r)) |$ is nondecreasing in $r$. It follows that, 
 $$ P_{\mathbb{H}}(F, B_{CC}(\hat{p}, r)) \leq 4r^{3}|B_{CC}(0, 1)| + \frac{2 \gamma}{\sqrt{h}} r^{4}|B_{CC}(0, 1)|$$
 and the claim follows using $r <\frac{C}{2\gamma} \sqrt{h} $.
  \end{proof}

 We are going now to recall  an $L^1$ estimate that holds  for sets which satisfy density estimates.
 
\begin{lemma}[Volume-distance comparison]\label{lemma-vol-dist}
Let $F\in \mathcal{P}$  bounded and  with the  following property: there exists   $\theta>0$ and $r_0<1$  such that for all $ p \in \partial (F(1))$ and $r \leq r_0$ there holds 
\begin{equation}\label{proprietaE} P_{\mathbb{H}}(F, B_{CC}(p, r) )=P_{\mathbb{H}}(F(1), B_{CC}(p, r) ) \geq \theta r^{3}.\end{equation}
Then for any $G \in \mathcal{P}$   and any $l>0$ we have
\begin{align*}
|F\triangle G| \leq  \omega_{CC} \frac{5^{4}}{\theta} l \max \left\{1,\frac{l^{3}}{r_0^3} \right\} P_{\mathbb{H}}(F)  +
\frac{2}{l} \int_{F \triangle G} d_{CC}(p, \partial F(1)) \,dp,
\end{align*}
where $\omega_{CC}:=|B_{CC}(p, 1)|$. 
\end{lemma}
\begin{remark}\label{osservazionoutward}\upshape
Note that due to Proposition \ref{density},   $E\in \mathcal{P}$  bounded, and Proposition \ref{propertyTh}(1) 
then  every  $F\in T_h(E)$  is bounded and satisfies \eqref{proprietaE} with $\theta$ depending only on universal constants and $r_0\leq C\sqrt{h}$. 

Therefore Lemma \ref{lemma-vol-dist} implies that   there exists  a constant $C>0$ depending only on the isoperimetric constant and on $\omega_{CC}$ 
such that for all $l\leq r_0\leq C \sqrt{h}$ and for all $F\in T_h(E)$ there holds 
\begin{equation}\label{stimal1} |F\triangle G| \leq  C\left( l P_{\mathbb{H}}(F)  +
\frac{1}{l} \int_{F \triangle G} d_{CC}(p, \partial F(1)) \,dp\right).
\end{equation} 
\end{remark}
\begin{proof}We follow \cite[Proposition~3.3]{BCK}. Define
\begin{align*}
V:=\left\{ p \in F \triangle G\,: \, d_{CC}(p, \partial F(1)) \leq \frac{l}{2} \right\}, \qquad W=\left\{p \in F \triangle G\,: \, d_{CC}(p, \partial F(1)) \geq \frac{l}{2} \right\}
\end{align*}
so that $|F \triangle G| \leq |V|+ |W| $. We have 
$$   \frac{2}{l} \int_{F \triangle G}  d_{CC}(p, \partial F(1)) \,dp\geq \frac{2}{l}  \int _{W}  d_{CC}(p, \partial F(1)) \,dp\geq |W|. $$
 Next,  by the compactness of $\partial F(1)$  and the definition of $V$,     there exists a finite  family of disjoint balls $B_{CC}(p_{j}, l)$ with  $p_{j}\in \partial F(1)$,  such that $V \subseteq \cup_{j=1}^mB_{CC}(p_{j}, 5l) $,  due to the covering theorem  in \cite[Theorem~1.2]{Heinonen}).
  If $l \geq r_0$ recalling the property \eqref{proprietaE} we get 
\begin{align*}
|V| &\leq \sum_{j=1}^{m} |B_{CC}(p_{j}, 5l)|= \sum_{j=1}^{m} \omega_{CC} (5l)^{4}=\omega_{CC} \frac{(5l)^{4}}{\theta r_0^3} \sum_{j=1}^{m} \theta r_0^3 \\& \leq \omega_{CC} \frac{(5l)^{4}}{\theta r_0^3} \sum_{j=1}^{m} P_{\mathbb{H}}(F, B_{CC}(p_{j},r_0) ) \leq \omega_{CC} \frac{(5l)^{4}}{\theta r_0^3}P_{\mathbb{H}}\big(F, \cup_{j=1}^{m} B_{CC}(p_{j},r_0) \big) \leq \omega_{CC} \frac{(5l)^{4}}{\theta r_0^3} P_{\mathbb{H}}(F)
\end{align*}
where we have used the fact that the balls $B_{CC}(p_{j}, r_0) \subset B_{CC}(p_{j}, l)$ are disjoint.
If $ l \leq r_0$  we obtain similarly
\begin{align*}
|V| &\leq \sum_{j=1}^{m} |B_{CC}(p_{j}, 5l)|=\omega_{CC} \frac{5^{4}l}{\theta} \sum_{j=1}^{m} \theta l^{3}\\& \leq \omega_{CC} \frac{5^{4}l}{\theta} \sum_{j=1}^{m} P_{\mathbb{H}}(F, B_{CC}(p_{j},l) ) \leq \omega_{CC} \frac{5^{4}l}{\theta}P_{\mathbb{H}}\big(F, \cup_{j=1}^{m} B_{CC}(p_{j},l) \big) \leq \omega_{CC} \frac{5^{4}l}{\theta} P_{\mathbb{H}}(F)
\end{align*}
and the claim follows.
\end{proof}
 
 \subsection{The flat flow} 
 We define the following discrete time scheme.  Take a  bounded set  $E\in \mathcal{P}$, and consider any $T_h(E)$  minimizer of \eqref{Jh}.
Then we construct the following piecewise constant flows of sets 
\[t\in [0,T]\longrightarrow E_h(t)\qquad E_h(t)=\begin{cases} E &t\in [0, h)\\
T_h(E)& t\in [h, 2h)\\ 
T_h(T_h(E))=:T^{2}_{h}(E)& t\in [2h, 3h)\\
\dots &\dots \end{cases} \]
In order to have a unique selection, we  choose each   $T^{k}(E)$  to be  $(T^{k}(E))^+$ (equivalently we could choose $T^{k}(E)$  to be  $(T^{k}(E))^-$).  
Note that by Proposition \ref{propertyTh}, item 4. the map 
 $t \rightarrow \Per(E_{h}(t))$ is monotone nonincreasing, that is 
 \begin{align}\label{monper}\Per(E_{h}(t)) \leq \Per (E_{h}(s))\qquad 0\leq s\leq t.\end{align}
Moreover $E_{h}(t)$ is bounded  for $t>0$ by  Proposition~\ref{propertyTh}(1) independently of $h$. 
 
 \begin{remark} \upshape Recalling Remark \ref{remarkdelta} we infer that \begin{align}\label{propertyJa}
\Per(T_{h}(E))  +\frac{1}{h}\int_{E \triangle (T_{h}(E))} dist(p, \partial E) dp\leq \Per(E).
\end{align}
By Lemma~\ref{lemma-vol-dist} and in particular by \eqref{stimal1} in Remark \ref{osservazionoutward}, , we obtain for  $l\leq C\sqrt{h}$, where $C$ is a constant depending only on the isoperimetric constant,   the following holds    for any $k\in \mathbb{N}$: \begin{align}\notag
|T^{k}_h(E) \triangle T^{k+1}_h(E)| &\leq  Cl P_{\mathbb{H}}(T^{k}(E)) +
\frac{C}{l} \int_{T^{k}(E) \triangle T^{k+1}(E)}  \mbox{dist}(p, \partial T^{k}_h(E)) \,dp\\
& \leq  C l  P_{\mathbb{H}}(T^{k}_h(E)) +
\frac{C}{l} h\Big( P_{\mathbb{H}}(T^{k}_h(E)) -P_{\mathbb{H}}(T^{k+1}(E)) \Big).\label{stimaiter}
 \end{align}
 
\end{remark}

\begin{remark}\label{rem112}\upshape
For sets $A,B,C$ there holds: $A \triangle C \subset (A \triangle B) \cup (B \triangle C) $.
\end{remark}

We provide the existence of the flat flow, as the time step $h\to 0$, starting from $E$  bounded.  
 
\begin{theorem}\label{thmflatflow} 
Let $E \in \mathcal{P}$ be a bounded set. 
There exists a sequence $h_{i} \to 0$, and a constant $C>0$ only depending on isoperimetric constant and on $p_0=\Per(E)$, such that the iteration scheme $t \to E_{h_{i}}(t)$ converges to a map $t \to E(t) \in \mathcal{P}$, $E(t)$ bounded and such that
\begin{align}\label{KVL1c}
|E(t) \triangle E(t')|&\leq C|t-t'|^{1/2}
\qquad \forall\,\, t,t' \in (0,\infty), \quad |t-t'|<1.
\end{align}
In addition, if $|\partial E|=0$,   then \eqref{KVL1c} holds also for  any $t,t' \in [0,\infty), |t-t'|<1$.
    \end{theorem}
\begin{proof}
Let    $0<t<t'$ such that $0<t'-t<1$ and let us fix $0<h\leq t'-t $ with $\lfloor\frac{t}{h}\rfloor\geq 1$,  where for $r\geq 0$ we denote $\lfloor r\rfloor:=\max\{n\in \N, n\leq r\}  $.   Let us introduce  $k, m \in 
\N$, such that  $k=\lfloor\frac{t}{h}\rfloor\geq 1$  and $k+m=\lfloor\frac{t'}{h}\rfloor\geq 1$.
We introduce the following notation $$[t]_{h}:=kh, \qquad [t']_{h}:=(k+m)h, \qquad h \leq mh=[t']_{h}-[t]_{h} 
\leq |t-t'|+h \leq 2|t-t'|. $$ Note that by definition $E_{h}(t)= E_{h}([t]_{h})=T^{k}_{h}(E)$, and $E_{h}(t')= E_{h}([t']_{h})=T^{k+m}_{h}(E)$.

We choose $l=\frac{Ch}{2|t-t'|^{1/2}}<C\sqrt{h}$, where $C$ is the constant in \eqref{stimaiter},  so that
$ml \leq C |t-t'|^{1/2}$. 
Now, we use estimate \eqref{stimaiter},  Remark \ref{rem112}, and   the fact that $ P_{\mathbb{H}}(T^{k}(E))\leq \Per(E)=p_0$  (see \eqref{monper}, to  get \begin{align}\label{KVL1} 
|E_{h}(t) \triangle E_{h}(t')|&=  |T^{k}_{h}(E) \triangle T^{k+m}_{h}(E)| \leq \sum_{j=1}^{m} |T^{k+j}_{h}(E) \triangle T^{k+j-1}_{h}(E)| \notag \\
& \leq Cl   p_0 m +\frac{C}{l}h \Big( P_{\mathbb{H}}(T^{k}(E)) -P_{\mathbb{H}}(T^{k+m}(E)) \Big) \notag\\
& \leq C^2 p_{0}|t-t'|^{1/2} + 2p_{0}|t-t'|^{1/2} =(C^2 p_0+2p_0)|t-t'|^{1/2}.
\end{align}  
By the a priori bounds given by  \eqref{monper}, and by the fact that $E_{h}(t)$ is bounded uniformly in $h$,  we may apply Theorem~\ref{CompactnessStrong}. to get that for all $t>0$, there exists a subsequence (depending on $t$)  $h_t\to 0$ such that 
 $E_h(t)$ converges in $L^1$ to some $E(t)$.
Using a simple diagonalization argument we can find a sequence $h_{i} \to 0$ such that
$$  \lim_{i \to \infty} E_{h_{i}}(t)=: E(t) \qquad \forall t \in \mathbb{Q} \cap (0,\infty).$$
From \eqref{KVL1} we infer that
\begin{align}\label{KVL1b}
|E(t) \triangle E(t')|&\leq (C^2p_0+2p_0) |t-t'|^{1/2}
\qquad \forall t,t' \in \mathbb{Q} \cap (0,\infty), |t-t'|<1.
\end{align}
Using again Theorem~\ref{CompactnessStrong} and \eqref{KVL1b} we extend univocally $E(t)$ to $(0,\infty)$ such that
\begin{align*}
|E(t) \triangle E(t')|&\leq  (C^2p_0+2p_0)  |t-t'|^{1/2}
\qquad \forall t,t' \in (0,\infty), |t-t'|<1
\end{align*}
holds, yielding \eqref{KVL1c}.\\
Finally let us consider the continuity of $E(t)$ at the initial time $t=0$. Obviously $E_{h}(0)=E$ for any $h >0$.
Notice that by \eqref{propertyJa} we have  
\begin{align*}
\int_{ E_{h_{i}}(h_{i}) \triangle E} dist (x, \partial E) dp \leq h_{i}\Per(E) \to 0 \text{ for } h_{i} \to 0.
\end{align*}
 On the other hand by  Theorem~\ref{CompactnessStrong} for a subsequence of $h_{i}$ (still denoted as $h_{i}$ for simplicity) we have that $\chi_{E_{h_{i}}(h_{i})} \to \chi_{B}$ in $L^{1} (\mathbb{H})$ for  some  $B \in \mathcal{P}$.  By Fatou's lemma we get 
 $$ \int_{   E\Delta  B} dist (p, \partial E) dp\leq \lim\inf_i\int_{E\Delta   E_{h_{i}}(h_{i})} dist (p, \partial E) dp  \leq  \lim\inf_i h_{i}\Per(E)  =0.$$
 This implies that $ dist(p, \partial E)=0$ for a.e. $p\in  B\setminus E$ and $p\in E\setminus B$, which implies that $B\subseteq \overline E$. Moreover  Since $|\partial E|=0$, it follows that $\chi_E=\chi_B$ that is   
    $\chi_{E\setminus E_{h_{i}}(h_{i}) }\to 0$ in $L^{1}(\mathbb{H})$.  Finally 
given $t \in (0,1)$ we compute
\begin{align*}
|E_{h_{i}}(t) \triangle E | \leq |E_{h_{i}}(t) \triangle E_{h_{i}}(h_{i}) |+ |E_{h_{i}}(h_{i}) \triangle E | \leq C|t-h_{i}|^{1/2} + |E\Delta E_{h_{i}}(h_{i})   | 
\end{align*}
and passing to the limit we obtain
$|E(t) \triangle E | \leq  C|t|^{1/2} $
as claimed.
\end{proof}
 

\subsection{Remarks on the consistency of the scheme} \label{sec:op}
Following \cite{capcitti2, ferrari, bc} a family of surfaces $t\in[0, T)\mapsto \Sigma_t=\partial E_t$, where     
$\Sigma_{t}$ is a  surface in $\mathbb{H}$,  given by the boundary of a $C^2$ set $E_t\subseteq \mathbb{H}$ for every $t\in [0, T)$, evolves by horizontal mean curvature flow if all its non-singular points $p(t) \in \Sigma_{t} \setminus (\Sigma_{t})_{0}$ move according to the velocity
\begin{equation} \label{geometricflow}\partial_{t} p(t) = - H_{\mathbb{H}}(p(t), E_{t}) \nu_{\mathbb{H}}(p(t))\end{equation} 
where $\nu_{\mathbb{H}}(p(t))$ is the horizontal exterior  normal to $E_t$ at the point $p(t)$.   
 If  $\Sigma_{t}$ can be described as the $0$-level set of a function $u$, that is $\Sigma_{t}=\{ (t,p) \in [0, T) \times \mathbb{H} \, : \, u(t,p)=0\}$, 
differentiation with respect to time of the equation  $0=u(t, p(t))$ yields (recalling \eqref{horgrad} and  \eqref{normal}) 
\begin{align*}
0 &=u_{t}(t, p(t)) +\nabla u(t, p(t)) \cdot \partial_{t} p(t)\\
&=u_{t}(t, p(t)) -\nabla u(t, p(t)) \cdot 
H_{\mathbb{H}}(p(t), E_{t}) \frac{\nabla_{\mathbb{H}}u (p(t))}{|\nabla_{\mathbb{H}}u (p(t))|_{g}}\\
&=u_{t} -  H_{\mathbb{H}} \nabla u \cdot 
 \frac{Xu\, X + Yu \,Y}{|\nabla_{\mathbb{H}}u|_{g}} = u_{t} -  H_{\mathbb{H}} \frac{ (Xu)^{2}+ (Yu)^{2}}{|\nabla_{\mathbb{H}}u|_{g}}\\
 & =u_{t}(t, p(t)) -H_{\mathbb{H}}(p(t), E_{t}) |\nabla_{\mathbb{H}}u(p(t))|_{g}.
\end{align*}
Therefore  away from the singular points 
$u$ is a solution to the PDE
\begin{equation}\label{levelsetpde} \frac{u_{t}}{|\nabla_{\mathbb{H}}u|_{g}} = H_{\mathbb{H}}(u).\end{equation}

 There are very few explicit  examples of smooth evolutions by horizontal mean curvature flow.  
 A first trivial observation is the following. Every   family of axially symmetric cylinders of the form    $C(r)_t=\{p=(x,y,z), \sqrt{x^2+y^2}=\sqrt{r^2-2t}\}$  is an evolution by horizontal mean curvature flow  for every $t\in [0, \frac{r^2}{2})$.  Indeed in every point of $C(r)_t=\partial E(r)_t=\partial\{p=(x,y,z), \sqrt{x^2+y^2}<\sqrt{r^2-2t}\}$  the horizontal mean curvature coincides with the euclidean mean curvature. 
 
 A less trivial example is the evolution of the Kor\'anyi sphere, obtained in \cite[Section 3.1]{ferrari}.  Recalling the definition given in Remark \ref{spherecurvature}, the Kor\'anyi sphere  is \[\partial B_K(0,R)=\{(x,y,z), \ (x^2+y^2)^2+16 z^2=R^4\} \]  and its 
 evolution by   horizontal mean curvature is 
\[\partial B^t_K(0,R)=\{(x,y,z),\  (x^2+y^2)^2+16 z^2+12 t(x^2+y^2)+12t^2=R^4\}.\]
for $t\in \left[0,  \frac{R^2}{2\sqrt{3}}\right)$. 
It is interesting to observe that despite the fact that $B_{K}(0,R) \in \mathcal{B}$ (see Definition~\ref{definitionb}), the sets $B^t_K(0,R)$  for $t\in \left(0,  \frac{R^2}{2\sqrt{3}}\right)$ do not belong to $\mathcal{B}$ since the horizontal mean curvature 
\[  H_{\mathbb{H}}(p,(B_K)_t(0,1))= \frac{3(x^2+y^2)+6t}{\sqrt{x^2+y^2}\sqrt{1+24t^2}}\] 
is positive but unbounded and $\delta_{\lambda} B_K^t(0, R)= B_K^{ \lambda^2 t} \left(0, \lambda R\right)$
holds.

  In general we cannot   expect to have smooth solutions to the horizontal  mean curvature flow, so weak notions of the evolutions have been introduced, mainly in the setting of level set formulation and viscosity solutions \cite{capcitti2, ferrari}. We recall also some results providing optimal control interpretation of the evolution by horizontal mean curvature flow, \cite{dirr, ddg, grande}.  In particular in  \cite{ferrari} existence and uniqueness of the viscosity solution to the level set equation associated to the horizontal mean curvature flow has been provided in the class of  axially symmetric sets. More recently   \cite{bc}   establishes a general uniqueness theorem and comparison principle for viscosity solutions to the horizontal  mean curvature flow, overcoming limitations near characteristic points by proving that continuous viscosity solutions can be derived as uniform limits of approximating Riemannian flows. 
  
  In this paper we provide the construction of the flat flow, Theorem \ref{thmflatflow}, which can be interpreted as another weak formulation of the horizontal mean curvature flow. 
  
 A natural question that we leave open is the consistency of the scheme we propose, that is the flat flow, with the smooth evolution by horizontal mean curvature flow and with other weak formulation of the flow, in particular with level set viscosity solutions. 
  
  The first problem would be the consistency with respect to smooth evolutions.  First of all, we observe that due to Theorem \ref{teo:7.1} 
  we may define the evolution by horizontal mean curvature flow for smooth sets as follows. 
   \begin{definition}\label{smoothhorizontal} We say that a family of surfaces $t\mapsto \Sigma_t=\partial E_t$, where  $E(t)$ are bounded $C^3$ sets  is a smooth horizontal mean curvature flow if
\[\partial_t d_{E(t)}(p)=\Div_{\mathbb{H}}\nabla_{\mathbb{H}} d_{E(t)}(p)\qquad \forall p\in \Sigma(t)\setminus \Sigma_0(t),\]
where $\Sigma_0(t)$ is the set of singular points in $\Sigma (t)$ and 
$d_E$ is the signed distance function from $E$ as defined in \eqref{sig} (with respect to the Carnot--Carathéodory distance $d_{CC}$).
\end{definition} 
We may define also smooth strict sub-and superflow of evolutions: 
 a family $\Sigma(t)=\partial E(t)$, where $E(t)$ denotes a $C^3$  bounded set in $\mathbb{H}$, is a strict superflow (resp. subflow) of the smooth horizontal mean curvature flow in $[t_0,t_1]$ if
 there exists $A\supseteq \cup_{t\in [t_0, t_1]}\Sigma(t)$, such that $d_{E(t)}$ is $C^3(A\setminus  \cup_{t\in [t_0, t_1]}\Sigma_0(t))$,  and  $\hat{\delta}>0$ for which 
\[\partial_t d_{E(t)}\geq \Div_{\mathbb{H}} \nabla_{\mathbb{H}} d_{E(t)}+\hat{\delta} \qquad \text{(resp. $\partial_t d_{E(t)}\leq \Div_{\mathbb{H}} \nabla_{\mathbb{H}} d_{E(t)}-\hat{\delta}$)} \] for $p\in A\setminus \cup_{t\in [t_0,t_1]} \Sigma_0(t)$, $t\in [t_0, t_1]$.

  We expect that if $\Sigma(t)=\partial E(t)$ is a smooth superflow, (resp. subflow)    then there exists $h_0$ such that for all $h\leq h_0$,  $t\in [t_0, t_1-h]$, it holds
\[E(t+h)\subseteq T_h(E(t))\qquad \text{ (resp. }E(t+h)\supseteq T_h(E(t))).\]
 
Such a result could be obtained in the case of evolving smooth surfaces without characteristic points, by adapting the calibration argument used in Proposition \ref{convexset}. Nonetheless the validity of a general result is at the moment unclear to us. We  recall also  that in \cite{cddg}  a new approximation and regularization of horizontal  mean curvature flow, similar to the Bence-Merriman-Osher algorithm, has been proposed,   and a consistency result  in that case was  obtained  for smooth  non-characteristic evolutions.

\renewcommand{\thesection}{}

\appendix\renewcommand{\thesection}{\Alph{section}}
\setcounter{equation}{0}
\renewcommand{\theequation}{\Alph{section}\arabic{equation}}

 
\end{document}